\documentclass[11pt]{article}
\usepackage{geometry}                
\usepackage{amsmath,amssymb}
\usepackage{amsthm}
\usepackage{hyperref}
\usepackage{graphicx,epsfig,color}
\usepackage{booktabs,bm,multirow,mathabx}
\usepackage{cases}
\usepackage[caption=false]{subfig}
\usepackage{lscape}
\newtheorem{theorem}{Theorem}
\newtheorem{proposition}[theorem]{Proposition}
\newtheorem{remark}[theorem]{Remark}

\def\mbbr{\mathbb R}
\def\l{\left} 
\def\r{\right} 
\def\mbf{\mathbf}
\def\rml{{\rm L}}
\def\rmc{{\rm C}}
\def\rmq{{\rm Q}}

\def\spa{{\rm span}}

\DeclareMathOperator*{\minimize}{minimize}

\def\wt{\widetilde}
\def\wh{\widehat}
\def\wc{\widecheck}

\newcommand{\beq}{\begin{equation}}\newcommand{\eeq}{\end{equation}}
\newcommand{\bit}{\begin{itemize}}\newcommand{\eit}{\end{itemize}}
\newcommand{\bem}{\begin{bmatrix}}\newcommand{\eem}{\end{bmatrix}}

\title{GPBiLQ and GPQMR: Two iterative methods for unsymmetric partitioned linear systems\thanks{This work was supported by the National Natural Science Foundation of China (No.12171403 and No.11771364), the Natural Science Foundation of Fujian Province of China (No.2020J01030), and the Fundamental Research Funds for the Central Universities (No.20720210032).}}

\author{Kui Du\thanks{Corresponding author. School of Mathematical Sciences and Fujian Provincial Key Laboratory of Mathematical Modeling and High Performance Scientific Computing, Xiamen University, Xiamen 361005, China (kuidu@xmu.edu.cn).},\quad Jia-Jun Fan\thanks{School of Mathematical Sciences, Xiamen University, Xiamen 361005, China ({fanjiajun@stu.xmu.edu.cn}).},\quad Fang Wang\thanks{School of Mathematical Sciences, Xiamen University, Xiamen 361005, China ({fangwang@stu.xmu.edu.cn}).}} 
\date{}                                           

\begin{document}
\maketitle

\begin{abstract} We introduce two iterative methods, GPBiLQ and GPQMR, for solving unsymmetric partitioned linear systems. The basic mechanism underlying GPBiLQ and GPQMR is a novel simultaneous tridiagonalization via biorthogonality that allows for short-recurrence iterative schemes. Similar to the biconjugate gradient method, it is possible to develop another method, GPBiCG, whose iterate (if it exists) can be obtained inexpensively from the GPBiLQ iterate. Whereas the iterate of GPBiCG may not exist, the iterates of GPBiLQ and GPQMR are always well defined as long as the biorthogonal tridiagonal reduction process does not break down. We discuss connections between the proposed methods and some existing methods, and give numerical experiments to illustrate the performance of the proposed methods.   

\vspace{2mm}  

\vspace{2mm}  

\end{abstract}

\section{Introduction} 
We are concerned with the partitioned linear system \beq\label{gplin}\bem\mbf M &\mbf A \\ \mbf B & \mbf N \eem\bem\mbf x\\ \mbf y\eem=\bem\mbf b\\ \mbf c\eem,\eeq where $\mbf M\in\mbbr^{m\times m}$, $\mbf N\in\mbbr^{n\times n}$, $\mbf A\in\mbbr^{m\times n}$, $\mbf B\in\mbbr^{n\times m}$, $\mbf b\in\mbbr^m$, and $\mbf c\in\mbbr^n$. System \eqref{gplin} arises in many areas (see, for example, \cite{elman1996fast,elman2014finit}). If $\mbf M$ and $\mbf N$ are invertible and the inversion operations are computationally inexpensive, then the equivalent right preconditioned system \beq\label{rgplin}\bem\mbf I &\mbf A\mbf N^{-1} \\ \mbf B\mbf M^{-1} & \mbf I \eem\bem\wt {\mbf x}\\ \wt {\mbf y}\eem=\bem\mbf M &\mbf A \\ \mbf B & \mbf N \eem\bem\mbf M^{-1} &\mbf 0 \\ \mbf 0 & \mbf N^{-1} \eem\bem\wt {\mbf x}\\ \wt {\mbf y}\eem=\bem\mbf b\\ \mbf c\eem,\quad \bem\mbf x\\ \mbf y\eem=\bem\mbf M^{-1} &\mbf 0 \\ \mbf 0 & \mbf N^{-1} \eem\bem\wt {\mbf x}\\ \wt {\mbf y}\eem\eeq can be solved instead of \eqref{gplin}. System \eqref{rgplin} is a special case of the following partitioned linear system \beq\label{lin} \bem\lambda \mbf I & \mbf A\\ \mbf B & \mu\mbf I \eem\bem\mbf x\\ \mbf y\eem=\bem\mbf b\\ \mbf c\eem.\eeq Note that $\lambda$ and/or $\mu$ may be zero. 

Recently Montoison and Orban \cite{montoison2023gpmr} introduced an iterative method named GPMR (general partitioned minimum residual) for solving \eqref{lin} with the assumption that $\mbf A$ and $\mbf B$ are nonzero, and $\mbf b$ and $\mbf c$ are nonzero. GPMR exploits the $2\times 2$ block structure of \eqref{lin} and is based on a simultaneous orthogonal Hessenberg reduction process. Mathematically, the approximation solution of GPMR is tantamount to the sum of the two approximation solutions of Block-GMRES applied to the following system with two right-hand sides \beq\label{blin}\bem\lambda \mbf I & \mbf A\\ \mbf B & \mu\mbf I \eem\bem\mbf x^b & \mbf x^c\\ \mbf y^b & \mbf y^c\eem=\bem\mbf b  & \mbf 0\\ \mbf 0 &\mbf c\eem.\eeq The storage and work per iteration of GPMR are similar to those of GMRES \cite{saad1986gmres}, and thus GPMR for \eqref{lin} outperforms Block-GMRES for \eqref{blin}. In the experiments given in \cite{montoison2023gpmr}, GPMR terminates significantly earlier than GMRES on a residual-based stopping criterion.    

The storage requirements of GMRES and GPMR increase with iterations. For large-scale problems, storage may become prohibitive after a few iterations. To overcome this issue, biorthogonalization methods (for example, BiLQ \cite{montoison2020bilq}, BiCG \cite{fletcher1976conju}, and QMR \cite{freund1991qmr}) based on three-term recurrences (the Lanczos biorthogonalization process \cite{lanczos1950itera}) can be used. The BiLQ and QMR iterates are always well defined at each iteration as long as the Lanczos biorthogonalization does not break down, but the BiCG iterate is not always well defined. It is possible to transition from the BiLQ iterate to the BiCG iterate whenever it is well defined. We emphasis that BiLQ and QMR are general iterative solvers which do not exploit the $2\times 2$ block structure of \eqref{lin}. In this paper, following the spirit of BiLQ and QMR, we develop two iterative methods named GPBiLQ and GPQMR for solving system \eqref{lin} by using the $2\times 2$ block structure and a simultaneous biorthogonal tridiagonal reduction process (see Section 2 for details).     

The paper is organized as follows. In the rest of this section, we give other related research and some notation. Section 2 gives background on the simultaneous orthogonal Hessenberg reduction process used in GPMR and the simultaneous biorthogonal tridiagonal reduction process used in GPBiLQ and GPQMR. Sections 3 and 4 discuss GPBiLQ and GPQMR for the partitioned linear system \eqref{lin}, respectively. It is possible to develop another method, GPBiCG, in the spirit of BiCG. We also give the relationship between the GPBiLQ iterate and the GPBiCG iterate in section 3. Numerical experiments are given in section 5, followed by extensions and concluding remarks in section 6.

{\it Related research}. There are a few solvers that have been tailored specifically to the block structure of systems \eqref{gplin}. Based on a simultaneous bidiagonalization procedure, Estrin and Greif \cite{estrin2018spmr} introduced SPMR, a family of saddle-point minimum residual solvers, for solving system \eqref{gplin} with $\mbf N=\mbf 0$. Based on Saunders--Simon--Yip tridiagonalization \cite{saunders1988two}, Buttari et al. \cite{buttari2019tridi} developed USYMLQR, a combination of the methods USYMLQ and USYMQR \cite{saunders1988two}, for solving system \eqref{gplin} with $\mbf B=\mbf A^\top$, $\mbf N=\mbf 0$, and symmetric positive definite $\mbf M$. When $\mbf B=\mbf A^\top$, $\mbf M$ is symmetric positive definite, and $\mbf N$ is symmetric negative definite, system \eqref{gplin} is said to be symmetric quasi-definite (SQD). Orban and Arioli \cite{orban2017itera} proposed two families of Krylov subspace methods for solving SQD linear systems. The first family is based on the generalized Golub--Kahan bidiagonalization process \cite{arioli2013gener}, including generalizations of well-known methods such as LSQR \cite{paige1982lsqr}, CRAIG \cite{craig1955nstep}, and LSMR \cite{fong2011lsmr} but also a new method named CRAIG-MR. The second family follows from a related Lanczos process with metric defined by blkdiag$(\mbf M, \mbf N)$ and contains methods operating directly on the SQD system, generalizing the CG \cite{hestenes1952metho} and MINRES \cite{paige1975solut} methods. Montoison and Orban \cite{montoison2021tricg} developed TriCG and TriMR, two methods related to Block-CG and Block-MINRES just like the relationship between GPMR and Block-GMRES, for solving SQD systems with $\mbf b$ and $\mbf c$ both nonzero. Compared with GPMR, the common feature of the above methods is that they all use short recurrences, so the storage and computational costs per iteration are fixed. Under some conditions (see Remark \ref{relationtossy} of section 2), the simultaneous biorthogonal tridiagonal reduction process used in GPBiLQ, GPBiCG, and GPQMR reduces to the Saunders--Simon--Yip tridiagonalization process, so GPBiCG coincides with TriCG, and GPQMR coincides with TriMR and GPMR in exact arithmetic. 

{\it Notation}. Lowercase (uppercase) boldface letters are reserved for column vectors (matrices) with the exception of $\mbf w_k$, which denotes a matrix with two columns. Lowercase lightface letters are reserved for scalars. For any vector $\mbf b\in\mbbr^n$, we use $\mbf b^\top$ and $\|\mbf b\|$ to denote the transpose and the Euclidean norm of $\mbf b$, respectively. We use $\mbf I_k$ to denote the $k\times k$ identity matrix, and use $\mbf e_k$ to denote the $k$th column of the identity matrix $\mbf I$ whose order is clear from the context. We use $\mbf 0$ to denote the zero vector (or matrix) of appropriate size. For any matrix $\mbf A\in\mbbr^{m\times n}$, we use $\mbf A^\top$ and $\|\mbf A\|$ to denote the transpose and the spectral norm of $\mbf A$, respectively. For nonsingular $\mbf A$, we use $\mbf A^{-1}$ to denote its inverse. Sometimes we use $\bigstar$ to denote a matrix whose elements are clear from the context.

\section{A simultaneous biorthogonal tridiagonal reduction process}

Montoison and Orban \cite{montoison2023gpmr} proved the following theorem. 
\begin{theorem}[{\cite[Theorem 2.1]{montoison2023gpmr}}]\label{SOHR} 
Let $\mbf A\in\mbbr^{m\times n}$ and $\mbf B\in\mbbr^{n\times m}$, and $p:=\min\{m,n\}$. There exist $\mbf V\in\mbbr^{m\times p}$ and $\mbf U\in\mbbr^{n\times p}$ with orthonormal columns, and upper Hessenberg $\mbf H\in\mbbr^{p\times p}$ and $\mbf F\in\mbbr^{p\times p}$ with nonnegative subdiagonal such that $$\mbf V^\top\mbf A\mbf U=\mbf H,\qquad \mbf U^\top\mbf B\mbf V=\mbf F.$$ 
\end{theorem}
Algorithm 1 is a transcription of Algorithm 2.1 of \cite{montoison2023gpmr}, which formalizes a simultaneous orthogonal Hessenberg reduction process derived from Theorem \ref{SOHR}. If $\mbf B=\mbf A^\top$, then Algorithm 1 reduces to the Saunders--Simon--Yip tridiagonalization process \cite{saunders1988two} in exact arithmetic. GPMR is based on Algorithm 1 and involves recurrences of growing length and growing cost.

\begin{table}[!htb]
\centering
\begin{tabular*}{170mm}{l}
\toprule {\bf Algorithm 1}: Simultaneous orthogonal Hessenberg reduction
\\ \hline\noalign{\smallskip} \quad {\bf Require}: $\mbf A$, $\mbf B$, $\mbf b$, $\mbf c$, all nonzero
\\ \noalign{\smallskip} \quad\hspace{.64mm} 1:\ \  $\beta\mbf v_1:=\mbf b$, $\gamma\mbf u_1:=\mbf c$ \hfill {\color[gray]{0.5}$\beta>0,\gamma>0$ so that $\|\mbf v_1\|=\|\mbf u_1\|=1$}  
\\ \noalign{\smallskip} \quad\hspace{.64mm} 2:\ \  {\bf for} $k=1,2,\cdots$  {\bf do}
\\ \noalign{\smallskip} \quad\hspace{.64mm} 3:\ \  \quad {\bf for} $i=1,2,\cdots, k $  {\bf do}
\\ \noalign{\smallskip}\hspace{.64mm} \quad 4:\ \  \quad\quad $h_{ik}=\mbf v_i^\top\mbf A\mbf u_k$
\\ \noalign{\smallskip}\hspace{.64mm} \quad 5:\ \  \quad\quad $f_{ik}=\mbf u_i^\top\mbf B\mbf v_k$
\\ \noalign{\smallskip} \quad\hspace{.64mm} 6:\ \  \quad {\bf end for}
\\ \noalign{\smallskip} \quad\hspace{.64mm} 7:\ \  \quad $h_{k+1,k}\mbf v_{k+1}=\mbf A\mbf u_k-\sum_{i=1}^kh_{ik}\mbf v_i$ \hfill  {\color[gray]{0.5}$h_{k+1,k}>0$ so that $\|\mbf v_{k+1}\|=1$}
\\ \noalign{\smallskip} \quad\hspace{.64mm} 8:\ \  \quad $f_{k+1,k}\mbf u_{k+1}=\mbf B\mbf v_k-\sum_{i=1}^kf_{ik}\mbf u_i$ \hfill  \hspace{48mm} {\color[gray]{0.5}$f_{k+1,k}>0$ so that $\|\mbf u_{k+1}\|=1$}
\\ \noalign{\smallskip} \quad\hspace{.64mm} 9:\ \  {\bf end for}\\ 
\bottomrule
\end{tabular*}
\end{table}

If we give up the use of orthogonal transformations in Algorithm 1, by biorthogonality and multiplications by $\mbf A^\top$, $\mbf B^\top$, as well as $\mbf A$ and $\mbf B$, we obtain a simultaneous biorthogonal tridiagonal reduction process, which is described as Algorithm 2. The relationship between Algorithm 1 and Algorithm 2 is analogous to that between the Arnoldi process \cite{arnoldi1951princ} and the Lanczos biorthogonalization process \cite{lanczos1950itera}.

\begin{table}[!htb]
\centering
\begin{tabular*}{170mm}{l}
\toprule {\bf Algorithm 2}: Simultaneous biorthogonal tridiagonal reduction
\\ \hline\noalign{\smallskip} \quad {\bf Require}: $\mbf A$, $\mbf B$, $\mbf b$, $\mbf c$, $\mbf f$, $\mbf g$, all nonzero
\\ \noalign{\smallskip} \quad\hspace{.64mm} 1:\ \ $\mbf p_0=\mbf 0$, $\mbf q_0=\mbf 0$, $\mbf u_0=\mbf 0$, $\mbf v_0=\mbf 0$
\\ \noalign{\smallskip} \quad\hspace{.64mm} 2:\ \  $\eta_1\mbf p_1=\wt{\mbf p}_1:=\mbf f$, $\beta_1\mbf q_1=\wt{\mbf q}_1:=\mbf b$ \hfill {\color[gray]{0.5}$(\eta_1,\beta_1)$ so that $\mbf p_1^\top\mbf q_1=1$}  
\\ \noalign{\smallskip} \quad\hspace{.64mm} 3:\ \  $\delta_1\mbf u_1=\wt{\mbf u}_1:=\mbf c$, $\gamma_1\mbf v_1=\wt{\mbf v}_1:=\mbf g$ \hfill  {\color[gray]{0.5}$(\delta_1,\gamma_1)$ so that $\mbf u_1^\top\mbf v_1=1$}  
\\ \noalign{\smallskip} \quad\hspace{.64mm} 4:\ \  {\bf for} $k=1,2,\cdots$  {\bf do}
\\ \noalign{\smallskip}\hspace{.64mm} \quad 5:\ \  \qquad $\alpha_k=\mbf p_k^\top\mbf A\mbf u_k$, $\theta_k=\mbf v_k^\top\mbf B\mbf q_k$
\\ \noalign{\smallskip} \quad\hspace{.64mm} 6:\ \  \qquad $\eta_{k+1}\mbf p_{k+1}=\wt{\mbf p}_{k+1}:=\mbf B^\top\mbf v_k-\delta_k\mbf p_{k-1}-\theta_k\mbf p_k$
\\ \noalign{\smallskip} \quad\hspace{.64mm} 7:\ \  \qquad $\beta_{k+1}\mbf q_{k+1}=\wt{\mbf q}_{k+1}:=\mbf A\mbf u_k-\gamma_k\mbf q_{k-1}-\alpha_k\mbf q_k$ \hfill  {\color[gray]{0.5}$(\eta_{k+1},\beta_{k+1})$ so that $\mbf p_{k+1}^\top\mbf q_{k+1}=1$}
\\ \noalign{\smallskip} \quad\hspace{.64mm} 8:\ \  \qquad $\delta_{k+1}\mbf u_{k+1}=\wt{\mbf u}_{k+1}:=\mbf B\mbf q_k-\eta_k\mbf u_{k-1}-\theta_k\mbf u_k$
\\ \noalign{\smallskip} \quad\hspace{.64mm} 9:\ \  \qquad $\gamma_{k+1}\mbf v_{k+1}=\wt{\mbf v}_{k+1}:=\mbf A^\top\mbf p_k-\beta_k\mbf v_{k-1}-\alpha_k\mbf v_k$ \hfill \hspace{15mm}  {\color[gray]{0.5}$(\delta_{k+1},\gamma_{k+1})$ so that $\mbf u_{k+1}^\top\mbf v_{k+1}=1$} 
\\ \noalign{\smallskip} \quad 10:\ \  {\bf end for}\\ 
\bottomrule
\end{tabular*}
\end{table}

For generic matrices $\mbf A$ and $\mbf B$, in exact arithmetic, Algorithm 2 will run to completion after $p=\min\{m,n\}$ steps, but for certain special matrices, as in other biorthogonalization methods (see, for example, \cite[Lecture 38]{trefethen1997numer}), lucky and unlucky breakdowns may happen. For simplicity, in this paper we do not consider breakdowns. The detailed discussion on all kinds of breakdowns and look-ahead techniques \cite{parlett1985look,freund1993imple} is left to future work. The scaling factors used in our implementation are $$\eta_k=|\wt{\mbf p}_k^\top\wt{\mbf q}_k|^{1/2},\quad\beta_k=\wt{\mbf p}_k^\top\wt{\mbf q}_k/\eta_k,\quad \delta_k=|\wt{\mbf u}_k^\top\wt{\mbf v}_k|^{1/2},\quad\gamma_k=\wt{\mbf u}_k^\top\wt{\mbf v}_k/\delta_k.$$ Define $$\mbf P_k:=\bem\mbf p_1& \mbf p_2& \cdots &\mbf p_k\eem,\quad \mbf Q_k:=\bem\mbf q_1& \mbf q_2& \cdots &\mbf q_k\eem,$$ $$\mbf U_k:=\bem\mbf u_1& \mbf u_2& \cdots &\mbf u_k\eem,\quad \mbf V_k:=\bem\mbf v_1& \mbf v_2& \cdots &\mbf v_k\eem,$$
 $$\mbf S_k:=\bem \alpha_1 & \gamma_2 & &\\ \beta_2 & \alpha_2 &\ddots & \\ & \ddots &\ddots & \gamma_k\\ &&\beta_k &\alpha_k \eem, \quad \mbf S_{k+1,k}:=\bem \mbf S_k\\ \beta_{k+1}\mbf e_k^\top\eem,\quad \mbf S_{k,k+1}:=\bem\mbf S_k & \gamma_{k+1}\mbf e_k\eem,$$ and $$\mbf T_k:=\bem \theta_1 & \eta_2 & &\\ \delta_2 & \theta_2 &\ddots & \\ & \ddots &\ddots & \eta_k\\ &&\delta_k &\theta_k \eem, \quad \mbf T_{k+1,k}:=\bem \mbf T_k\\ \delta_{k+1}\mbf e_k^\top\eem,\quad \mbf T_{k,k+1}:=\bem\mbf T_k & \eta_{k+1}\mbf e_k\eem.$$ After $k$ iterations of Algorithm 2, we have the following relations: 
\begin{subequations}\label{matrixrelation}\begin{align}
\mbf A \mbf U_k&=\mbf Q_{k+1}\mbf S_{k+1,k}=\mbf Q_k\mbf S_k+\beta_{k+1}\mbf q_{k+1}\mbf e_k^\top,
\\ \mbf A^\top\mbf P_k &=\mbf V_{k+1}\mbf S_{k,k+1}^\top=\mbf V_k\mbf S_k^\top+\gamma_{k+1}\mbf v_{k+1}\mbf e_k^\top,
\\ \mbf B \mbf Q_k &=\mbf U_{k+1}\mbf T_{k+1,k}=\mbf U_k\mbf T_k+\delta_{k+1}\mbf u_{k+1}\mbf e_k^\top,
\\ \mbf B^\top\mbf V_k &=\mbf P_{k+1}\mbf T_{k,k+1}^\top=\mbf P_k\mbf T_k^\top+\eta_{k+1}\mbf p_{k+1}\mbf e_k^\top,
\\ \mbf P_k^\top\mbf Q_k& =\mbf U_k^\top\mbf V_k=\mbf I_k.
\end{align}\end{subequations}
For $i\geq 1$, we can show by induction that 
\begin{align*}
 \mbf p_{2i} & \in\spa\{\mbf f, \ldots, (\mbf B^\top\mbf A^\top)^{i-1}\mbf f, \mbf B^\top\mbf g, \ldots, (\mbf B^\top\mbf A^\top)^{i-1}\mbf B^\top\mbf g\}, 
 \\ \mbf p_{2i+1} & \in\spa\{\mbf f, \ldots, (\mbf B^\top\mbf A^\top)^i\mbf f, \mbf B^\top\mbf g, \ldots, (\mbf B^\top\mbf A^\top)^{i-1}\mbf B^\top\mbf g\},
 \\	\mbf q_{2i} & \in\spa\{\mbf b, \ldots, (\mbf A\mbf B)^{i-1}\mbf b, \mbf A\mbf c, \ldots, (\mbf A\mbf B)^{i-1}\mbf A\mbf c\}, 
 \\ \mbf q_{2i+1} & \in\spa\{\mbf b, \ldots, (\mbf A\mbf B)^i\mbf b, \mbf A\mbf c, \ldots, (\mbf A\mbf B)^{i-1}\mbf A\mbf c\},
 \\	\mbf u_{2i} & \in\spa\{\mbf c, \ldots, (\mbf B\mbf A)^{i-1}\mbf c, \mbf B\mbf b, \ldots, (\mbf B\mbf A)^{i-1}\mbf B\mbf b\}, 
 \\ \mbf u_{2i+1} & \in\spa\{\mbf c, \ldots, (\mbf B\mbf A)^i\mbf c, \mbf B\mbf b, \ldots, (\mbf B\mbf A)^{i-1}\mbf B\mbf b\},
 \\ \mbf v_{2i} & \in\spa\{\mbf g, \ldots, (\mbf A^\top\mbf B^\top)^{i-1}\mbf g, \mbf A^\top\mbf f, \ldots, (\mbf A^\top\mbf B^\top)^{i-1}\mbf A^\top\mbf f\}, 
 \\ \mbf v_{2i+1} & \in\spa\{\mbf g, \ldots, (\mbf A^\top\mbf B^\top)^i\mbf g, \mbf A^\top\mbf f, \ldots, (\mbf A^\top\mbf B^\top)^{i-1}\mbf A^\top\mbf f\}.
\end{align*} This means the subspaces generated by Algorithm 2 can be viewed as the union of four block Krylov subspaces generated by $\mbf B^\top\mbf A^\top$, $\mbf A\mbf B$, $\mbf B\mbf A$, and $\mbf A^\top\mbf B^\top$ with respective starting blocks $\bem\mbf f &\mbf B^\top\mbf g\eem$, $\bem\mbf b &\mbf A\mbf c \eem$, $\bem \mbf c&\mbf B\mbf b\eem$, and $\bem \mbf g & \mbf A^\top\mbf f\eem$. When $\mbf B=\mbf A^\top$, $\mbf f = \mbf b$, and $\mbf g=\mbf c$, we have the following remark. 

\begin{remark}\label{relationtossy}
If $\mbf B=\mbf A^\top$, $\mbf f = \mbf b$, and $\mbf g=\mbf c$, then for each $k\geq 1$ we have $\mbf p_k=\mbf q_k$ and $\mbf u_k=\mbf v_k$ in exact arithmetic. Thus Algorithm 2 reduces to the Saunders--Simon--Yip tridiagonalization process \cite{saunders1988two}, which has the following relations 	\begin{align*}
\mbf A \mbf U_k&=\mbf Q_{k+1}\mbf S_{k+1,k},
\quad \mbf A^\top\mbf Q_k =\mbf U_{k+1}\mbf S_{k,k+1}^\top,\quad \mbf S_k=\mbf Q_k^\top\mbf A\mbf U_k,\quad \mbf Q_k^\top\mbf Q_k =\mbf U_k^\top\mbf U_k=\mbf I_k.
\end{align*}
\end{remark}

Now we consider projecting the coefficient matrix of system \eqref{lin} onto a smaller subspace where the projected matrix has a block structure with tridiagonal off-diagonal blocks. By \eqref{matrixrelation}, we have $$\bem\mbf 0 &\mbf A\\ \mbf B &\mbf 0\eem\bem \mbf Q_k & \mbf 0\\ \mbf 0& \mbf U_k\eem=\bem \mbf Q_{k+1} & \mbf 0\\ \mbf 0& \mbf U_{k+1}\eem\bem \mbf 0 & \mbf S_{k+1,k}\\ \mbf T_{k+1,k} &  \mbf 0\eem,$$ 
which implies \beq\label{proj}\bem\lambda\mbf I & \mbf A\\ \mbf B & \mu\mbf I\eem\bem \mbf Q_k & \mbf 0\\ \mbf 0& \mbf U_k\eem=\bem \mbf Q_k & \mbf 0\\ \mbf 0& \mbf U_k \eem \bem \lambda\mbf I & \mbf S_k\\ \mbf T_k & \mu\mbf I \eem + \bem \mbf q_{k+1} & \mbf 0\\ \mbf 0 & \mbf u_{k+1} \eem \bem \mbf 0 & \beta_{k+1}\mbf e_k^\top\\ \delta_{k+1}\mbf e_k^\top &  \mbf 0\eem.\eeq
Let $$ \mbf \Pi_k:=\bem \mbf e_1& \mbf e_{k+1} & \cdots & \mbf e_i & \mbf e_{i+k} & \cdots & \mbf e_k & \mbf e_{2k}\eem$$ denote the permutation introduced by Paige \cite{paige1974bidia}. Define $$\mbf W_k:=\bem \mbf Q_k & \mbf 0\\ \mbf 0& \mbf U_k\eem\mbf \Pi_k= \bem\mbf w_1 & \mbf w_2 &\cdots & \mbf w_k\eem, \quad \mbf w_k:=
\bem\mbf q_k & \mbf 0\\  \mbf 0 & \mbf u_k\eem.$$  
Multiplying both sides of equation \eqref{proj} by $\mbf \Pi_k$ and using $\mbf \Pi_k\mbf \Pi_k^\top=\mbf I$, we obtain \begin{align}\label{proj2}\bem\lambda\mbf I & \mbf A\\ \mbf B & \mu\mbf I\eem\mbf W_k =\mbf W_k\mbf \Pi_k^\top\bem \lambda\mbf I & \mbf S_k\\ \mbf T_k & \mu\mbf I \eem\mbf \Pi_k+ \bem \mbf q_{k+1} & \mbf 0\\ \mbf 0 & \mbf u_{k+1} \eem \bem \mbf 0 & \beta_{k+1}\mbf e_k^\top\\ \delta_{k+1}\mbf e_k^\top &  \mbf 0\eem\mbf \Pi_k=\mbf W_{k+1}\mbf H_{k+1,k},\end{align} 
where $$\mbf H_{k+1,k}=\bem \mbf E_{11} & \mbf E_{12} &  & \\ \mbf E_{21} & \mbf E_{22} & \ddots & \\ & \ddots & \ddots & \mbf E_{k-1,k} \\ & & \mbf E_{k,k-1} & \mbf E_{kk}\\ &&& \mbf E_{k+1,k} \eem,\ \ \mbf E_{i,i-1}=\bem 0 & \beta_i \\ \delta_i & 0\eem, \ \mbf E_{ii}=\bem \lambda & \alpha_i \\ \theta_i & \mu\eem, \ \mbf E_{i-1,i}=\bem 0 & \gamma_i \\ \eta_i & 0\eem.$$
Using equation \eqref{proj2} and the principles similar in spirit to BiLQ and QMR, we can derive two iterative methods, named GPBiLQ and GPQMR, for solving system \eqref{lin}. They are the subjects of the following two sections.

\section{Derivations of GPBiLQ and GPBiCG} 
Define $$\mbf H_{k-1,k}:=\bem \mbf E_{11} & \mbf E_{12} &  & & \\ \mbf E_{21} & \mbf E_{22} & \ddots & & \\ & \ddots & \ddots & \mbf E_{k-2,k-1} & \\ & & \mbf E_{k-1,k-2} & \mbf E_{k-1,k-1} &\mbf E_{k-1,k} \eem\in\mbbr^{(2k-2)\times 2k}.$$ We define GPBiLQ as the method that generates an approximation \beq\label{bilqsol}\bem\mbf x_k^\rml\\ \mbf y_k^\rml\eem=\mbf W_k\mbf z_k^\rml,\eeq 
where $\mbf z_k^\rml\in\mbbr^{2k}$ solves \beq\label{bilq}\minimize_{\mbf z\in\mbbr^{2k}}\|\mbf z\|\quad \mbox{ subject to }\quad \mbf H_{k-1,k}\mbf z=\beta_1\mbf e_1+\delta_1\mbf e_2.\eeq  For the matrix $\mbf H_{k-1,k}$, we have the following proposition. 

\begin{proposition}\label{prop1} For each $k$, the matrix $\mbf H_{k-1,k}$ has full row rank as long as Algorithm 2 does not break down.
\end{proposition}
\proof If Algorithm 2 does not break down, then for $i=2,3\ldots,k$, all $\gamma_i$ and $\eta_i$ are nonzero. Hence, all $\mbf E_{i-1,i}$ are nonsingular. It follows that $$\bem \mbf E_{12} &  & & \\  \mbf E_{22} & \mbf E_{23} & & \\ \mbf E_{32} & \ddots & \ddots& & \\ &\ddots &\ddots & \mbf E_{k-2,k-1} & \\ & & \mbf E_{k-1,k-2} & \mbf E_{k-1,k-1} &\mbf E_{k-1,k} \eem\in\mbbr^{(2k-2)\times(2k-2)}$$ is nonsingular. That is to say, $\mbf H_{k-1,k}\in\mbbr^{(2k-2)\times 2k}$ has a submatrix of rank $2k-2$, and thus $\mbf H_{k-1,k}$ has full row rank. 
\endproof

Proposition \ref{prop1} implies that the constraint of (\ref{bilq}) is consistent as long as Algorithm 2 does not break down, and thus the GPBiLQ iterate is always well defined. 

We compute $\mbf z_k^\rml$ via the LQ factorization of $\mbf H_{k-1,k}$, which can be obtained from the LQ factorization \beq\label{lqHk}\mbf H_k=\bem \mbf E_{11} & \mbf E_{12} &  & \\ \mbf E_{21} & \mbf E_{22} & \ddots & \\ & \ddots & \ddots & \mbf E_{k-1,k} \\ & & \mbf E_{k,k-1} & \mbf E_{k,k}\eem=\wt{\mbf L}_k\wt{\mbf Q}_k\eeq with \beq\label{lklk1}\wt{\mbf L}_k=\bem\hspace{-1mm}\begin{array}{cccccc|cc} \rho_1 & &&&&&&  \\ \nu_2 & \ddots &&&&&& \\ \omega_3 & \ddots &\ddots &&&&& \\ \zeta_4 &\ddots &\ddots &\ddots &&&& \\ \xi_5 &\ddots &\ddots &\ddots &\ddots &&& \\ &\ddots &\ddots &\ddots &\ddots &\rho_{2k-2} && \\ \hline &&\xi_{2k-1} &\zeta_{2k-1} &\omega_{2k-1} &\nu_{2k-1} & \ddot\rho_{2k-1} &\\ &&&\xi_{2k} &\zeta_{2k} &\omega_{2k} &\ddot\nu_{2k} & \ddot\rho_{2k}\end{array}\hspace{-1mm}\eem=\bem \mbf L_{k-1} & \mbf 0\\ \bigstar &  {\bigstar} \eem \in\mbbr^{2k\times 2k}\eeq and $$\wt{\mbf Q}_k^\top=\mbf G_{1,4}\mbf G_{3,6}\cdots\mbf G_{2k-5,2k-2}\mbf G_{2k-3,2k}\wt{\mbf G}_{2k-1,2k}\in\mbbr^{2k\times 2k}.$$ 
For $i=1,2,\ldots,k-1$, the structure of $\mbf G_{2i-1,2i+2}$ is $$\mbf G_{2i-1,2i+2}=\bem\mbf I_{2i-2}&&&&&\\& \times &\times &\times &\times &\\&\times &\times &\times &\times &\\&\times &\times &\times &\times &\\ &\times &\times &\times &\times & \\ &&&&& \mbf I_{2k-2i-2}\eem.$$ The $4\times 4$ diagonal block is the following product of four Givens rotations: $$\bem c_{i,1} &&& -s_{i,1}\\ & 1&&\\ &&1&\\ s_{i,1} & & & c_{i,1}\eem \bem c_{i,2}&-s_{i,2}&& \\ s_{i,2}&c_{i,2}&&\\ &&1&\\ &&&1 \eem\bem 1&&& \\ &c_{i,3}&&-s_{i,3}\\ &&1&\\ &s_{i,3}&&c_{i,3} \eem \bem 1&&&\\ &c_{i,4}& -s_{i,4}& \\ & s_{i,4}& c_{i,4}&\\ &&&1\eem.$$ The structure of $\wt{\mbf G}_{2k-1,2k}$ is \beq\label{strug}\wt{\mbf G}_{2k-1,2k}=\bem\mbf I_{2k-2} &&\\ &c_k&-s_k\\ &s_k& c_k\eem.\eeq 
Indeed, equation \eqref{lqHk} yields the LQ factorization \beq\label{lqHk1}\mbf H_{k-1,k}=\bem\mbf I_{2(k-1)} &\mbf 0\eem\mbf H_k=\bem\mbf I_{2(k-1)} &\mbf 0\eem\wt {\mbf L}_k\wt{\mbf Q}_k=\bem\mbf L_{k-1} & \mbf 0\eem\wt{\mbf Q}_k=\bem\mbf L_{k-1} & \mbf 0\eem\mbf G_k^\top,\eeq where $$\mbf G_k:=\mbf G_{1,4}\mbf G_{3,6}\cdots\mbf G_{2k-5,2k-2}\mbf G_{2k-3,2k}.$$ 
Note that $$\mbf L_k=\bem \rho_1 & &&&&&& \\ \nu_2 & \ddots &&&&&& \\  \omega_3 & \ddots &\ddots &&&&& \\ \zeta_4 &\ddots &\ddots &\ddots &&&& \\ \xi_5 &\ddots &\ddots &\ddots &\ddots &&& \\ &\ddots &\ddots &\ddots &\ddots &\ddots && \\ &&\xi_{2k-1} &\zeta_{2k-1} &\omega_{2k-1} &\nu_{2k-1} & \rho_{2k-1} &\\ &&&\xi_{2k} &\zeta_{2k} &\omega_{2k} &\nu_{2k} & \rho_{2k}\eem\in\mbbr^{2k\times 2k}$$ differs from $\wt{\mbf L}_k$ only in the $(2k-1,2k-1)$th, $(2k,2k-1)$th, and $(2k,2k)$th elements. The computations for the elements of the matrices $\wt{\mbf L}_k$, $\mbf G_{2i-1,2i+2}$, and $\wt{\mbf G}_{2k-1,2k}$ are given in Appendix \ref{lq}. If $$\mbf t_{k-1}:=\bem \varpi_1 & \varpi_2 & \cdots  & \varpi_{2k-3} & \varpi_{2k-2}\eem^\top$$ solves \beq\label{lk1}\mbf L_{k-1}\mbf t=\beta_1\mbf e_1+\delta_1\mbf e_2,\eeq then by \eqref{lqHk1} the solution $\mbf z_k^\rml$ of \eqref{bilq} is $$\mbf z_k^\rml=\mbf G_k\bem\mbf t_{k-1}\\ \mbf 0\eem.$$ Using the structure of $\mbf L_{k-1}$, we get $$\mbf t_{k-1}=\bem\mbf t_{k-2}^\top & \varpi_{2k-3} & \varpi_{2k-2}\eem^\top.$$ The computations for the elements of $\mbf t_{k-1}$ are given in Appendix \ref{wtt}.  

The $k$th GPBiLQ iterate is $\mbf W_k\mbf z_k^\rml$. To avoid storing $\mbf W_k$, we define $$\mbf F_k:=\mbf W_k\mbf G_k=\bem \mbf f_1^x & \cdots & \mbf f_{2k-3}^x& \mbf f_{2k-2}^x & \wt{\mbf f}_{2k-1}^x & \wt{\mbf f}_{2k}^x\\ \mbf f_1^y & \cdots & \mbf f_{2k-3}^y& \mbf f_{2k-2}^y & \wt{\mbf f}_{2k-1}^y & \wt{\mbf f}_{2k}^y\eem.$$ 
Using $$\mbf W_k=\bem \mbf W_{k-1} & \mbf w_k\eem, \quad \mbf G_k=\bem \mbf G_{k-1} & \\ & \mbf I_2\eem\mbf G_{2k-3,2k},$$ we get \beq\label{wg}\mbf W_k\mbf G_k=\bem\mbf W_{k-1}\mbf G_{k-1} &\mbf w_k\eem\mbf G_{2k-3,2k}.\eeq 
Recall that $\mbf G_{2k-3,2k}\in\mbbr^{2k\times 2k}$ is a block diagonal matrix with the first diagonal block being $\mbf I_{2k-4}$. This means that the first $2k-4$ columns of $\mbf F_k$ and $\mbf F_{k-1}$ are identical. Thus the $k$th GPBiLQ iterate satisfies \begin{align*}\bem\mbf x_k^\rml\\ \mbf y_k^\rml\eem & =\mbf W_k\mbf G_k\mbf G_k^\top\mbf z_k^\rml=\mbf F_k\bem\mbf t_{k-1}\\ \mbf 0\eem=\bem \mbf f_1^x & \cdots & \mbf f_{2k-3}^x& \mbf f_{2k-2}^x & \wt{\mbf f}_{2k-1}^x & \wt{\mbf f}_{2k}^x\\ \mbf f_1^y & \cdots & \mbf f_{2k-3}^y& \mbf f_{2k-2}^y & \wt{\mbf f}_{2k-1}^y & \wt{\mbf f}_{2k}^y\eem\bem\mbf t_{k-2}\\ \varpi_{2k-3} \\ \varpi_{2k-2} \\ 0 \\ 0\eem\\ &= \mbf F_{k-1}\bem\mbf t_{k-2}\\ \mbf 0\eem +\bem \mbf f_{2k-3}^x& \mbf f_{2k-2}^x \\ \mbf f_{2k-3}^y& \mbf f_{2k-2}^y \eem \bem \varpi_{2k-3} \\ \varpi_{2k-2} \eem= \bem\mbf x_{k-1}^\rml\\ \mbf y_{k-1}^\rml\eem+\bem \mbf f_{2k-3}^x& \mbf f_{2k-2}^x \\ \mbf f_{2k-3}^y& \mbf f_{2k-2}^y \eem \bem \varpi_{2k-3} \\ \varpi_{2k-2} \eem. \end{align*} 
The update is initialized with $\mbf x_1^\rml=\mbf 0$ and $\mbf y_1^\rml=\mbf 0$. Equating the last four columns on both sides of \eqref{wg}, we get \beq\label{eq:f4}\bem \mbf f_{2k-3}^x& \mbf f_{2k-2}^x & \wt{\mbf f}_{2k-1}^x & \wt{\mbf f}_{2k}^x\\ \mbf f_{2k-3}^y& \mbf f_{2k-2}^y & \wt{\mbf f}_{2k-1}^y & \wt{\mbf f}_{2k}^y\eem=\bem \wt{\mbf f}_{2k-3}^x& \wt{\mbf f}_{2k-2}^x & \mbf q_k & \mbf 0\\  \wt{\mbf f}_{2k-3}^y& \wt{\mbf f}_{2k-2}^y & \mbf 0 & \mbf u_k\eem\bem\mbf 0 & \mbf I_4\eem\mbf G_{2k-3,2k}\bem\mbf 0\\ \mbf I_4\eem.\eeq This yields the recursion for $\mbf f_{2k-3}^x$, $\mbf f_{2k-2}^x$, $\mbf f_{2k-3}^y$, and $\mbf f_{2k-2}^y$.

\subsection{From GPBiLQ to GPBiCG}
We define GPBiCG as the method that generates an approximation \beq\label{bicgsol}\bem\mbf x_k^\rmc\\ \mbf y_k^\rmc\eem=\mbf W_k\mbf z_k^\rmc,\eeq where $\mbf z_k^\rmc\in\mbbr^{2k}$ solves \beq\label{bicg}\mbf H_k\mbf z=\beta_1\mbf e_1+\delta_1\mbf e_2 .\eeq Just like BiCG, the above system maybe inconsistent, so that the GPBiCG iterate is not always well defined. 

Assume that $\wt{\mbf L}_k$ in the LQ factorization \eqref{lqHk} is nonsingular. If $\wt{\mbf t}_k$ solves $$\wt{\mbf L}_k\mbf t=\beta_1\mbf e_1+\delta_1\mbf e_2,$$ then it follows from \eqref{lklk1} and the fact that $\mbf t_{k-1}$ solves \eqref{lk1} that $\wt{\mbf t}_k=\bem \mbf t_{k-1}^\top & \wt\varpi_{2k-1} & \wt\varpi_{2k}\eem^\top$. The computations for $\wt\varpi_{2k-1}$ and $\wt\varpi_{2k}$ are given in Appendix \ref{wtt}. Using the LQ factorization $\mbf H_k=\wt{\mbf L}_k\wt{\mbf Q}_k$, we get $$\mbf z_k^\rmc=\wt{\mbf Q}_k^\top\wt{\mbf t}_k.$$ 
Define $\wt{\mbf F}_k:=\mbf W_k\wt{\mbf Q}_k^\top$. Using $\wt{\mbf Q}_k^\top=\mbf G_k\wt{\mbf G}_{2k-1,2k}$ and $\mbf F_k=\mbf W_k\mbf G_k$, we get $$\wt{\mbf F}_k=\mbf F_k \wt{\mbf G}_{2k-1,2k}=\bem \mbf f_1^x & \cdots & \mbf f_{2k-3}^x& \mbf f_{2k-2}^x & \wt{\mbf f}_{2k-1}^x & \wt{\mbf f}_{2k}^x\\ \mbf f_1^y & \cdots & \mbf f_{2k-3}^y& \mbf f_{2k-2}^y & \wt{\mbf f}_{2k-1}^y & \wt{\mbf f}_{2k}^y\eem\wt{\mbf G}_{2k-1,2k}.$$ 
The structure of $\wt{\mbf G}_{2k-1,2k}$ \eqref{strug} yields that only the last two columns of $\wt{\mbf F}_k$ and $\mbf F_k$ are different. Thus the $k$th GPBiCG iterate is \begin{align*}\bem\mbf x_k^\rmc\\ \mbf y_k^\rmc\eem  &=\mbf W_k\mbf z_k^\rmc=\mbf W_k\wt{\mbf Q}_k^\top \wt{\mbf t}_k=\wt{\mbf F}_k\wt{\mbf t}_k= \bem\mbf x_k^\rml\\ \mbf y_k^\rml\eem+\bem  \wt{\mbf f}_{2k-1}^x & \wt{\mbf f}_{2k}^x\\ \wt{\mbf f}_{2k-1}^y & \wt{\mbf f}_{2k}^y\eem\bem c_k & -s_k\\ s_k & c_k\eem\bem \wt\varpi_{2k-1}\\ \wt\varpi_{2k}\eem.\end{align*} This means that it is possible to transfer from the $k$th GPBiLQ iterate to the $k$th GPBiCG iterate cheaply provided that $\wt{\mbf L}_k$ is nonsingular. 

Algorithm 3 summarizes GPBiLQ, including the GPBiCG  iterate. 
  
\subsection{Residual estimates} 

The $k$th residual of GPBiLQ is $$\mbf r_k^\rml =\bem\mbf b\\ \mbf c\eem-\bem\lambda\mbf I & \mbf A\\ \mbf B & \mu\mbf I\eem\bem \mbf x_k^\rml\\ \mbf y_k^\rml\eem.$$ 
Using equations \eqref{proj2} and \eqref{bilqsol}, $\beta_1\mbf q_1=\mbf b$, and $\delta_1\mbf u_1=\mbf c$, we obtain \begin{align*}\mbf r_k^\rml =\mbf W_{k+1}(\beta_1\mbf e_1+\delta_1\mbf e_2)-\bem\lambda\mbf I & \mbf A\\ \mbf B & \mu\mbf I\eem\mbf W_k\mbf z_k^\rml = \mbf W_{k+1}(\beta_1\mbf e_1+\delta_1\mbf e_2)-\mbf W_{k+1}\mbf H_{k+1,k}\mbf z_k^\rml.\end{align*} 
Substituting $\mbf H_{k-1,k}\mbf z_k^\rml=\beta_1\mbf e_1+\delta_1\mbf e_2$ into the above equation yields \begin{align*} \mbf r_k^\rml= -\bem \mbf q_k  &\mbf 0 \\ \mbf 0 & \mbf u_k \eem \l(\mbf E_{k,k-1}\bem\mbf e_{2k-3}^\top\\ \mbf e_{2k-2}^\top \eem + \mbf E_{k,k}\bem\mbf e_{2k-1}^\top\\ \mbf e_{2k}^\top \eem\r) \mbf z_k^\rml- \bem \mbf q_{k+1}  &\mbf 0 \\ \mbf 0 & \mbf u_{k+1} \eem \mbf E_{k+1,k}\bem\mbf e_{2k-1}^\top\\ \mbf e_{2k}^\top \eem\mbf z_k^\rml.\end{align*}
Define $$\bem \vartheta_k \\ \varrho_k \eem := \l(\mbf E_{k,k-1}\bem\mbf e_{2k-3}^\top\\ \mbf e_{2k-2}^\top \eem + \mbf E_{k,k}\bem\mbf e_{2k-1}^\top\\ \mbf e_{2k}^\top \eem\r) \mbf z_k^\rml,$$ and $$\bem\chi_k\\ \varsigma_k\eem:= \mbf E_{k+1,k}\bem\mbf e_{2k-1}^\top\\ \mbf e_{2k}^\top \eem\mbf z_k^\rml.$$ We have $$\|\mbf r_k^\rml\|=\l\|\bem \vartheta_k\mbf q_k+\chi_k\mbf q_{k+1}\\ \varrho_k\mbf u_k +\varsigma_k\mbf u_{k+1}\eem\r\|.$$
The computations for $\vartheta_k$, $\varrho_k$, $\chi_k$, and $\varsigma_k$ are given in Appendix \ref{trcs}. 

The $k$th residual of GPBiCG is $$\mbf r_k^\rmc = \bem\mbf b\\ \mbf c\eem-\bem\lambda\mbf I & \mbf A\\ \mbf B & \mu\mbf I\eem\bem \mbf x_k^\rmc\\ \mbf y_k^\rmc\eem.$$ 
Using equations \eqref{proj2} and \eqref{bicgsol}, $\beta_1\mbf q_1=\mbf b$, and $\delta_1\mbf u_1=\mbf c$, we obtain \begin{align*}
\mbf r_k^\rmc &=\mbf W_{k+1}(\beta_1\mbf e_1+\delta_1\mbf e_2)-\bem\lambda\mbf I & \mbf A\\ \mbf B & \mu\mbf I\eem\mbf W_k\mbf z_k^\rmc = \mbf W_{k+1}(\beta_1\mbf e_1+\delta_1\mbf e_2)-\mbf W_{k+1}\mbf H_{k+1,k}\mbf z_k^\rmc.	
\end{align*}
Substituting $\mbf H_k\mbf z_k^\rmc=\beta_1\mbf e_1+\delta_1\mbf e_2$ into the above equation yields
$$\mbf r_k^\rmc = - \bem \mbf q_{k+1}  &\mbf 0 \\ \mbf 0 & \mbf u_{k+1} \eem \mbf E_{k+1,k}\bem\mbf e_{2k-1}^\top\\ \mbf e_{2k}^\top \eem\mbf z_k^\rmc.$$
Define $$\bem\wt\chi_k\\ \wt\varsigma_k\eem:= \mbf E_{k+1,k}\bem\mbf e_{2k-1}^\top\\ \mbf e_{2k}^\top \eem\mbf z_k^\rmc.$$ We have $$\|\mbf r_k^\rmc\|=\sqrt{|\wt\chi_k|^2\|\mbf q_{k+1}\|^2+|\wt\varsigma_k|^2\|\mbf u_{k+1}\|^2}.$$ The computations for $\wt\chi_k$ and $\wt\varsigma_k$ are given in Appendix \ref{tcs}.

\subsection{Computational expense and storage} The computational expense per iteration of GPBiLQ is dominated by the four matrix-vector products ($\mbf A\mbf u$, $\mbf A^\top\mbf p$, $\mbf B\mbf q$, $\mbf B^\top\mbf v$) and some vector operations. Note that $\mbf A$, $\mbf A^\top$, $\mbf B$, and $\mbf B^\top$ do not need to be formed explicitly and we only require four procedures to compute the four matrix-vector products. 

If in-place {\tt gemv} updates (for example, $\mbf p\leftarrow\mbf B^\top\mbf v-\delta\mbf p$) are available, GPBiLQ requires nine $m$-vectors ($\mbf x_k^\rml$, $\mbf p_{k-1}$, $\mbf p_k$, $\mbf q_{k-1}$, $\mbf q_k$, $ \mbf f_{2k-3}^x$,  $\mbf f_{2k-2}^x$, $\wt{\mbf f}_{2k-1}^x$, $\wt{\mbf f}_{2k}^x$) and nine $n$-vectors ($\mbf y_k^\rml$, $\mbf u_{k-1}$, $\mbf u_k$, $\mbf v_{k-1}$, $\mbf v_k$, $ \mbf f_{2k-3}^y$,  $\mbf f_{2k-2}^y$, $\wt{\mbf f}_{2k-1}^y$, $\wt{\mbf f}_{2k}^y$). If in-place {\tt gemv} updates are not available, additional $m$- and $n$-vectors are required to store $\mbf A\mbf u$, $\mbf B^\top\mbf v$, $\mbf A^\top\mbf p$, and $\mbf B\mbf q$. We remark that the storage requirement of GPBiLQ is comparable to that of GPMR(9) (a restarted variant of GPMR) since each iteration of GPMR requires one $m$-vector and one $n$-vector (see Table 3.1 of \cite{montoison2023gpmr}).

\begin{table}[!htb]
\centering
\begin{tabular*}{170mm}{l}
\toprule {\bf Algorithm 3}: GPBiLQ
\\ \hline\noalign{\smallskip} \quad {\bf Require}: $\mbf A$, $\mbf B$, $\mbf b$, $\mbf c$, $\mbf f$, $\mbf g$, all nonzero; $\lambda$, $\mu$, \tt maxit
\\ \noalign{\smallskip} \quad\hspace{.64mm} 1:\ \  $\eta_1\mbf p_1=\wt{\mbf p}_1:=\mbf f$, $\beta_1\mbf q_1=\wt{\mbf q}_1:=\mbf b$ \hfill {\color[gray]{0.5}$(\eta_1,\beta_1)$ so that $\mbf p_1^\top\mbf q_1=1$}
\\ \noalign{\smallskip} \quad\hspace{.64mm} 2:\ \  $\delta_1\mbf u_1=\wt{\mbf u}_1:=\mbf c$, $\gamma_1\mbf v_1=\wt{\mbf v}_1:=\mbf g$ \hfill {\color[gray]{0.5}$(\delta_1,\gamma_1)$ so that $\mbf u_1^\top\mbf v_1=1$}  
\\ \noalign{\smallskip} \quad\hspace{.64mm} 3:\ \ $\alpha_1=\mbf p_1^\top\mbf A\mbf u_1$, $\theta_1=\mbf v_1^\top\mbf B\mbf q_1$ 
\\ \noalign{\smallskip} \quad\hspace{.64mm} 4:\ \ $\eta_2\mbf p_2=\wt{\mbf p}_2:=\mbf B^\top\mbf v_1-\theta_1\mbf p_1$, $\beta_2\mbf q_2=\wt{\mbf q}_2:=\mbf A\mbf u_1-\alpha_1\mbf q_1$ \hfill {\color[gray]{0.5}$(\eta_2,\beta_2)$ so that $\mbf p_2^\top\mbf q_2=1$}   
\\ \noalign{\smallskip}\hspace{.64mm} \quad 5:\ \ $\delta_2\mbf u_2=\wt{\mbf u}_2:=\mbf B\mbf q_1-\theta_1\mbf u_1$, $\gamma_2\mbf v_2=\wt{\mbf v}_2:=\mbf A^\top\mbf p_1-\alpha_1\mbf v_1$ \hfill {\color[gray]{0.5}$(\delta_2,\gamma_2)$ so that $\mbf u_2^\top\mbf v_2=1$}
\\ \noalign{\smallskip} \quad\hspace{.64mm} 6:\ \ $\mbf x_1^\rml=\mbf 0$, $\mbf y_1^\rml=\mbf 0$, $\wt{\mbf f}_1^x=\mbf q_1$, $\wt{\mbf f}_2^x=\mbf 0$, $\wt{\mbf f}_1^y=\mbf 0$, $\wt{\mbf f}_2^y=\mbf u_1$
\\ \noalign{\smallskip} \quad\hspace{.64mm} 7:\ \ $\bar\rho_1=\lambda$, $\bar\alpha_1=\alpha_1$, $\bar\nu_2=\theta_1$, $\bar\rho_2=\mu$, $\bar\omega_3=0$, $\bar\nu_3=\beta_2$, $\bar\zeta_4=\delta_2$
\\ \noalign{\smallskip} \quad\hspace{.64mm} 8:\ \ $k=2$
\\ \noalign{\smallskip} \quad\hspace{.64mm} 9:\ \ {\bf while} $k\leq\tt maxit$ {\bf do} 
\\ \noalign{\smallskip} \quad 10:\ \  \quad $\alpha_k=\mbf p_k^\top\mbf A\mbf u_k$, $\theta_k=\mbf v_k^\top\mbf B\mbf q_k$
\\ \noalign{\smallskip} \quad 11:\ \  \quad $\eta_{k+1}\mbf p_{k+1}=\wt{\mbf p}_{k+1}:=\mbf B^\top\mbf v_k-\delta_k\mbf p_{k-1}-\theta_k\mbf p_k$
\\ \noalign{\smallskip} \quad 12:\ \  \quad $\beta_{k+1}\mbf q_{k+1}=\wt{\mbf q}_{k+1}:=\mbf A\mbf u_k-\gamma_k\mbf q_{k-1}-\alpha_k\mbf q_k$ \hfill  {\color[gray]{0.5}$(\eta_{k+1},\beta_{k+1})$ so that $\mbf p_{k+1}^\top\mbf q_{k+1}=1$}
\\ \noalign{\smallskip} \quad 13:\ \  \quad $\delta_{k+1}\mbf u_{k+1}=\wt{\mbf u}_{k+1}:=\mbf B\mbf q_k-\eta_k\mbf u_{k-1}-\theta_k\mbf u_k$
\\ \noalign{\smallskip} \quad 14:\ \  \quad $\gamma_{k+1}\mbf v_{k+1}=\wt{\mbf v}_{k+1}:=\mbf A^\top\mbf p_k-\beta_k\mbf v_{k-1}-\alpha_k\mbf v_k$ \hfill  {\color[gray]{0.5}$(\delta_{k+1},\gamma_{k+1})$ so that $\mbf u_{k+1}^\top\mbf v_{k+1}=1$} 
\\ \noalign{\smallskip} \quad 15:\ \  \quad Compute $\rho_{2k-3}$, $\nu_{2k-2}$, $\rho_{2k-2}$, $\omega_{2k-1}$, $\nu_{2k-1}$, $\zeta_{2k}$, $\omega_{2k}$, $\xi_{2k+1}$, $\zeta_{2k+1}$, $\xi_{2k+2}$ \hfill {\color[gray]{0.5} (see Appendix \ref{lq})}
\\ \noalign{\smallskip} \quad 16:\ \  \quad Compute $\bar\rho_{2k-1}$, $\bar\alpha_k$, $\bar\nu_{2k}$, $\bar\rho_{2k}$, $\bar\omega_{2k+1}$, $\bar\nu_{2k+1}$, $\bar\zeta_{2k+2}$  \hfill {\color[gray]{0.5} (see Appendix \ref{lq})}
\\ \noalign{\smallskip} \quad 17:\ \  \quad Compute $\varpi_{2k-3}$, $\varpi_{2k-2}$ \hfill {\color[gray]{0.5} (see Appendix \ref{wtt})}
\\ \noalign{\smallskip} \quad 18:\ \  \quad Compute $ \mbf f_{2k-3}^x$,  $\mbf f_{2k-2}^x$, $\wt{\mbf f}_{2k-1}^x$, $\wt{\mbf f}_{2k}^x$, $\mbf f_{2k-3}^y$, $\mbf f_{2k-2}^y$, $\wt{\mbf f}_{2k-1}^y$, $\wt{\mbf f}_{2k}^y$ \hfill {\color[gray]{0.5} (use \eqref{eq:f4})}
\\ \noalign{\smallskip} \quad 19:\ \  \quad $\mbf x_k^\rml=\mbf x_{k-1}^\rml+\varpi_{2k-3}\mbf f_{2k-3}^x+ \varpi_{2k-2}\mbf f_{2k-2}^x$ 
\\ \noalign{\smallskip} \quad 20:\ \  \quad $\mbf y_k^\rml=\mbf y_{k-1}^\rml+ \varpi_{2k-3} \mbf f_{2k-3}^y +\varpi_{2k-2}\mbf f_{2k-2}^y$
\\ \noalign{\smallskip} \quad 21:\ \  \quad {\bf if} $\bar\rho_{2k-1}\bar\rho_{2k}-\bar\alpha_k\bar\nu_{2k}\neq 0$ {\bf then}
\\ \noalign{\smallskip} \quad 22:\ \  \qquad Compute $c_k$, $s_k$, $\wt\varpi_{2k-1}$, $\wt\varpi_{2k}$\hfill {\color[gray]{0.5} (see Appendix \ref{lq} and Appendix \ref{wtt})} 
\\ \noalign{\smallskip} \quad 23:\ \  \qquad $\mbf x_k^\rmc=\mbf x_k^\rml + (c_k\wt\varpi_{2k-1} -s_k\wt\varpi_{2k})\wt{\mbf f}_{2k-1}^x +(s_k\wt\varpi_{2k-1}+c_k\wt\varpi_{2k}) \wt{\mbf f}_{2k}^x $
\\ \noalign{\smallskip} \quad 24:\ \  \qquad $\mbf y_k^\rmc=\mbf y_k^\rml + (c_k\wt\varpi_{2k-1} -s_k\wt\varpi_{2k})\wt{\mbf f}_{2k-1}^y +(s_k\wt\varpi_{2k-1}+c_k\wt\varpi_{2k}) \wt{\mbf f}_{2k}^y $
\\ \noalign{\smallskip} \quad 25:\ \  \quad  {\bf end if}
\\ \noalign{\smallskip} \quad 26:\ \  \quad $k=k+1$
\\ \noalign{\smallskip} \quad 27:\ \  {\bf end while}\\ 
\bottomrule
\end{tabular*}
\end{table}

\section{Derivation of GPQMR}

We define GPQMR as the method that generates an approximation \beq\label{qmrsol}\bem\mbf x_k^\rmq\\ \mbf y_k^\rmq\eem=\mbf W_k\mbf z_k^\rmq,\eeq where $\mbf z_k^\rmq\in\mbbr^{2k}$ solves \beq\label{qmr}\minimize_{\mbf z\in\mbbr^{2k}}\|\mbf H_{k+1,k}\mbf z-(\beta_1\mbf e_1+\delta_1\mbf e_2)\|.\eeq For the matrix $\mbf H_{k+1,k}$, we have the following proposition. 

\begin{proposition}\label{prop2} For each $k$, the matrix $\mbf H_{k+1,k}$ has full column rank as long as Algorithm 2 does not break down.
\end{proposition}
\proof If Algorithm 2 does not break down, then for $i=2,3\ldots,k+1$, all $\beta_i$ and $\delta_i$ are nonzero. Hence, all $\mbf E_{i,i-1}$ are nonsingular. It follows that $$\bem \mbf E_{21} & \mbf E_{22} & \mbf E_{23} & & \\  & \mbf E_{32} & \ddots & \ddots & \\  && \ddots & \ddots & \mbf E_{k-1,k} \\ && & \mbf E_{k,k-1} & \mbf E_{kk}\\ & &&& \mbf E_{k+1,k} \eem \in\mbbr^{2k\times 2k}$$ is nonsingular. That is to say, $\mbf H_{k+1,k}\in\mbbr^{(2k+2)\times 2k}$ has a submatrix of rank $2k$, and thus $\mbf H_{k+1,k}$ has full column rank. 
\endproof

Proposition \ref{prop2} implies that the least squares problem \eqref{qmr} has a unique solution as long as Algorithm 2 does not break down, and thus the GPQMR iterate is always well defined.

We compute $\mbf z_k^\rmq$ via the QR factorization \beq\label{qrhk1}\mbf H_{k+1,k}=\wh{\mbf Q}_k\bem\wh{\mbf R}_k\\ \mbf 0\eem,\eeq where $$\wh{\mbf R}_k=\bem \rho_1 & \nu_1 & \omega_1 & \zeta_1 & \xi_1 &&\\ & \ddots & \ddots & \ddots & \ddots & \ddots &\\ && \ddots & \ddots & \ddots & \ddots & \xi_{2k-4} \\ &&& \ddots & \ddots & \ddots & \zeta_{2k-3}\\ &&&& \ddots & \ddots & \omega_{2k-2}\\ &&&&& \ddots & \nu_{2k-1}\\ &&&&&& \rho_{2k} \eem\in\mbbr^{2k\times 2k}$$ and $$\wh{\mbf Q}_k^\top=\wh{\mbf G}_{2k-1,2k+2}\cdots\wh{\mbf G}_{3,6}\wh{\mbf G}_{1,4}\in\mbbr^{(2k+2)\times(2k+2)}.$$
For $i=1,2,\ldots,k$, the structure of $\wh{\mbf G}_{2i-1,2i+2}$ is $$\wh{\mbf G}_{2i-1,2i+2}=\bem\mbf I_{2i-2}&&&&&\\& \times &\times &\times &\times &\\&\times &\times &\times &\times &\\&\times &\times &\times &\times &\\ &\times &\times &\times &\times & \\ &&&&& \mbf I_{2k-2i}\eem.$$ 
The $4\times 4$ diagonal block is the following product of four Givens rotations: $$\bem 1&&&\\ &c_{i,4}& s_{i,4}& \\ & -s_{i,4}& c_{i,4}&\\ &&&1 \eem  \bem 1&&& \\ &c_{i,3}&&s_{i,3}\\ &&1&\\ &-s_{i,3}&&c_{i,3} \eem \bem c_{i,2}&s_{i,2}&& \\ -s_{i,2}&c_{i,2}&&\\ &&1&\\ &&&1 \eem
\bem c_{i,1} &&& s_{i,1}\\ & 1&&\\ &&1&\\ -s_{i,1} & & & c_{i,1}\eem.$$
The computations for the elements of the matrices $\wh{\mbf R}_k$ and $\wh{\mbf G}_{2i-1,2i+2}$ are given in Appendix \ref{qr}. We use $\varpi_1, \ldots, \varpi_{2k} $, $\bar\varpi_{2k+1}$ and $\bar\varpi_{2k+2}$ to represent the $2k+2$ elements of the vector $\wh{\mbf Q}_k^\top(\beta_1\mbf e_1+\delta_1\mbf e_2)$, i.e., $$\wh{\mbf Q}_k^\top(\beta_1\mbf e_1+\delta_1\mbf e_2)=\bem \varpi_1 &\varpi_2& \ldots &\varpi_{2k-1}&\varpi_{2k} & \bar\varpi_{2k+1} & \bar\varpi_{2k+2}\eem^\top.$$ The computations for the elements of $\wh{\mbf Q}_k^\top(\beta_1\mbf e_1+\delta_1\mbf e_2)$ are given in Appendix \ref{bart}. Define $$\mbf t_k:=\bem \varpi_1 &\varpi_2& \ldots &\varpi_{2k-1}&\varpi_{2k} \eem^\top.$$ 
Using the structure of $\wh{\mbf Q}_k^\top$, we get $$\mbf t_k=\bem\mbf t_{k-1}^\top & \varpi_{2k-1}&\varpi_{2k}\eem^\top.$$ 
By the QR factorization \eqref{qrhk1}, the solution of \eqref{qmr} satisfies \beq\label{rzt}\wh{\mbf R}_k\mbf z_k^\rmq=\mbf t_k.\eeq

The $k$th GPQMR iterate is $\mbf W_k\mbf z_k^\rmq$. To avoid storing $\mbf W_k$, we define $$\mbf F_k:=\mbf W_k\wh{\mbf R}_k^{-1}=\bem \mbf f_1^x & \cdots & \mbf f_{2k-3}^x& \mbf f_{2k-2}^x & \mbf f_{2k-1}^x & \mbf f_{2k}^x\\ \mbf f_1^y & \cdots & \mbf f_{2k-3}^y& \mbf f_{2k-2}^y & \mbf f_{2k-1}^y & \mbf f_{2k}^y\eem.$$ 
Equating the last two columns on both sides of $\mbf W_k=\mbf F_k\wh{\mbf R}_k$, we get \begin{align*}
\mbf f_{2k-1}^x &= (\mbf q_k -\xi_{2k-5}\mbf f_{2k-5}^x-\zeta_{2k-4}\mbf f_{2k-4}^x-\omega_{2k-3}\mbf f_{2k-3}^x-\nu_{2k-2}\mbf f_{2k-2}^x)/\rho_{2k-1},\\
\mbf f_{2k-1}^y &= (-\xi_{2k-5}\mbf f_{2k-5}^y-\zeta_{2k-4}\mbf f_{2k-4}^y-\omega_{2k-3}\mbf f_{2k-3}^y-\nu_{2k-2}\mbf f_{2k-2}^y)/\rho_{2k-1},\\ 
\mbf f_{2k}^x &= (-\xi_{2k-4}\mbf f_{2k-4}^x-\zeta_{2k-3}\mbf f_{2k-3}^x-\omega_{2k-2}\mbf f_{2k-2}^x-\nu_{2k-1}\mbf f_{2k-1}^x)/\rho_{2k},\\
\mbf f_{2k}^y &= (\mbf u_k -\xi_{2k-4}\mbf f_{2k-4}^y-\zeta_{2k-3}\mbf f_{2k-3}^y-\omega_{2k-2}\mbf f_{2k-2}^y-\nu_{2k-1}\mbf f_{2k-1}^y)/\rho_{2k}.
\end{align*} 
Using $$\mbf W_k=\bem\mbf W_{k-1} & \mbf w_k\eem,\quad \wh{\mbf R}_k=\bem\wh{\mbf R}_{k-1} & \bigstar \\ & \bigstar\eem,$$ we get $$\mbf F_k=\bem \mbf F_{k-1} & \bigstar \eem.$$
Thus the $k$th GPQMR iterate satisfies \begin{align*}\bem\mbf x_k^\rmq\\ \mbf y_k^\rmq\eem & =\mbf W_k\wh{\mbf R}_k^{-1}\wh{\mbf R}_k\mbf z_k^\rmq=\mbf F_k\mbf t_k=\bem \mbf f_1^x & \cdots & \mbf f_{2k-3}^x& \mbf f_{2k-2}^x & \mbf f_{2k-1}^x & \mbf f_{2k}^x\\ \mbf f_1^y & \cdots & \mbf f_{2k-3}^y& \mbf f_{2k-2}^y & \mbf f_{2k-1}^y & \mbf f_{2k}^y\eem\mbf t_k \\ & = \bem\mbf x_{k-1}^\rmq\\ \mbf y_{k-1}^\rmq\eem+\bem \mbf f_{2k-1}^x& \mbf f_{2k}^x \\ \mbf f_{2k-1}^y& \mbf f_{2k}^y \eem \bem \varpi_{2k-1} \\ \varpi_{2k} \eem. \end{align*}
The update is initialized with $\mbf x_0^\rmq=\mbf 0$ and $\mbf y_0^\rmq=\mbf 0$.

Algorithm 4 summarizes GPQMR.

\subsection{Residual estimates}
The $k$th residual of GPQMR is $$\mbf r_k^\rmq =\bem\mbf b\\ \mbf c\eem-\bem\lambda\mbf I & \mbf A\\ \mbf B & \mu\mbf I\eem\bem \mbf x_k^\rmq\\ \mbf y_k^\rmq\eem.$$
Using equations \eqref{proj2} and \eqref{qmrsol}, $\beta_1\mbf q_1=\mbf b$, and $\delta_1\mbf u_1=\mbf c$, we obtain $$\mbf r_k^\rmq =\mbf W_{k+1}(\beta_1\mbf e_1+\delta_1\mbf e_2)-\bem\lambda\mbf I & \mbf A\\ \mbf B & \mu\mbf I\eem\mbf W_k\mbf z_k^\rmq = \mbf W_{k+1}(\beta_1\mbf e_1+\delta_1\mbf e_2)-\mbf W_{k+1}\mbf H_{k+1,k}\mbf z_k^\rmq.$$
Using equations \eqref{qrhk1} and \eqref{rzt}, and $\wh{\mbf Q}_k^\top(\beta_1\mbf e_1+\delta_1\mbf e_2)=\bem\mbf t_k^\top & \bar\varpi_{2k+1} & \bar\varpi_{2k+2} \eem^\top$, we obtain \begin{align*} \mbf r_k^\rmq  = \mbf W_{k+1}\wh{\mbf Q}_k\l(\wh{\mbf Q}_k^\top(\beta_1\mbf e_1+\delta_1\mbf e_2)-\bem\wh{\mbf R}_k\\ \mbf 0\eem\mbf z_k^\rmq\r)=\mbf W_{k+1}\wh{\mbf Q}_k( \bar\varpi_{2k+1}\mbf e_{2k+1} + \bar\varpi_{2k+2}\mbf e_{2k+2}).\end{align*} 
Then we have $$\|\mbf r_k^\rmq\|\leq\|\mbf W_{k+1}\|\sqrt{\bar\varpi_{2k+1}^2+\bar\varpi_{2k+2}^2}.$$

\subsection{Computational expense and storage} The computational expense and storage of GPQMR are almost the same as those of GPBiLQ. 

\begin{table}[!htb]
\centering
\begin{tabular*}{170mm}{l}
\toprule {\bf Algorithm 4}: GPQMR
\\ \hline\noalign{\smallskip} \quad {\bf Require}: $\mbf A$, $\mbf B$, $\mbf b$, $\mbf c$, $\mbf f$, $\mbf g$, all nonzero; $\lambda$, $\mu$,  \tt maxit
\\ \noalign{\smallskip} \quad\hspace{.64mm} 1:\ \ $\mbf p_0=\mbf 0$, $\mbf q_0=\mbf 0$, $\mbf u_0=\mbf 0$, $\mbf v_0=\mbf 0$
\\ \noalign{\smallskip} \quad\hspace{.64mm} 2:\ \  $\eta_1\mbf p_1=\wt{\mbf p}_1:=\mbf f$, $\beta_1\mbf q_1=\wt{\mbf q}_1:=\mbf b$ \hfill {\color[gray]{0.5}$(\eta_1,\beta_1)$ so that $\mbf p_1^\top\mbf q_1=1$}
\\ \noalign{\smallskip} \quad\hspace{.64mm} 3:\ \  $\delta_1\mbf u_1=\wt{\mbf u}_1:=\mbf c$, $\gamma_1\mbf v_1=\wt{\mbf v}_1:=\mbf g$ \hfill {\color[gray]{0.5}$(\delta_1,\gamma_1)$ so that $\mbf u_1^\top\mbf v_1=1$}  
\\ \noalign{\smallskip} \quad\hspace{.64mm} 4:\ \ $\mbf x_0^\rmq=\mbf 0$, $\mbf y_0^\rmq=\mbf 0$, $\mbf f_{-3}^x=\mbf 0$, $\mbf f_{-2}^x=\mbf 0$, $\mbf f_{-1}^x=\mbf 0$, $\mbf f_0^x=\mbf 0$, $\mbf f_{-3}^y=\mbf 0$, $\mbf f_{-2}^y=\mbf 0$, $\mbf f_{-1}^y=\mbf 0$, $\mbf f_0^y=\mbf 0$
\\ \noalign{\smallskip} \quad\hspace{.64mm} 5:\ \ $k=1$
\\ \noalign{\smallskip} \quad\hspace{.64mm} 6:\ \ {\bf while} $k\leq\tt maxit$ {\bf do} 
\\ \noalign{\smallskip} \quad \hspace{.64mm} 7:\ \  \quad $\alpha_k=\mbf p_k^\top\mbf A\mbf u_k$, $\theta_k=\mbf v_k^\top\mbf B\mbf q_k$
\\ \noalign{\smallskip} \quad \hspace{.64mm} 8:\ \  \quad $\eta_{k+1}\mbf p_{k+1}=\wt{\mbf p}_{k+1}:=\mbf B^\top\mbf v_k-\delta_k\mbf p_{k-1}-\theta_k\mbf p_k$
\\ \noalign{\smallskip} \quad \hspace{.64mm} 9:\ \  \quad $\beta_{k+1}\mbf q_{k+1}=\wt{\mbf q}_{k+1}:=\mbf A\mbf u_k-\gamma_k\mbf q_{k-1}-\alpha_k\mbf q_k$ \hfill  {\color[gray]{0.5}$(\eta_{k+1},\beta_{k+1})$ so that $\mbf p_{k+1}^\top\mbf q_{k+1}=1$}
\\ \noalign{\smallskip} \quad 10:\ \  \quad $\delta_{k+1}\mbf u_{k+1}=\wt{\mbf u}_{k+1}:=\mbf B\mbf q_k-\eta_k\mbf u_{k-1}-\theta_k\mbf u_k$
\\ \noalign{\smallskip} \quad 11:\ \  \quad $\gamma_{k+1}\mbf v_{k+1}=\wt{\mbf v}_{k+1}:=\mbf A^\top\mbf p_k-\beta_k\mbf v_{k-1}-\alpha_k\mbf v_k$ \hfill  {\color[gray]{0.5}$(\delta_{k+1},\gamma_{k+1})$ so that $\mbf u_{k+1}^\top\mbf v_{k+1}=1$} 
\\ \noalign{\smallskip} \quad 12:\ \ \quad {\bf if} $k=1$ {\bf then}
\\ \noalign{\smallskip} \quad 13:\ \ \qquad $\bar\rho_1=\lambda$, $\bar\theta_1=\theta_1$, $\bar\nu_1=\alpha_1$, $\bar\rho_2=\mu$, $\bar\omega_1=0$, $\bar\nu_2=\eta_2$, $\bar\zeta_1=\gamma_2$, $\bar\varpi_1=\beta_1$, $\bar\varpi_2=\delta_1$
\\ \noalign{\smallskip} \quad 14:\ \ \quad {\bf end if}
\\ \noalign{\smallskip} \quad 15:\ \  \quad Compute $\rho_{2k-1}$, $\nu_{2k-1}$, $\rho_{2k}$, $\omega_{2k-1}$, $\nu_{2k}$, $\zeta_{2k-1}$, $\omega_{2k}$, $\xi_{2k-1}$, $\zeta_{2k}$, $\xi_{2k}$ \hfill \hspace{8.7mm} {\color[gray]{0.5} (see Appendix \ref{qr})}
\\ \noalign{\smallskip} \quad 16:\ \  \quad Compute $\bar\rho_{2k+1}$, $\bar\theta_{k+1}$, $\bar\nu_{2k+1}$, $\bar\rho_{2k+2}$, $\bar\omega_{2k+1}$, $\bar\nu_{2k+2}$, $\bar\zeta_{2k+1}$ \hfill {\color[gray]{0.5} (see Appendix \ref{qr})}
\\ \noalign{\smallskip} \quad 17:\ \  \quad Compute $\varpi_{2k-1}$, $\varpi_{2k}$, $\bar\varpi_{2k+1}$, $\bar\varpi_{2k+2}$ \hfill {\color[gray]{0.5} (see Appendix \ref{bart})}
\\ \noalign{\smallskip} \quad 18:\ \  \quad $\mbf f_{2k-1}^x= (\mbf q_k -\xi_{2k-5}\mbf f_{2k-5}^x-\zeta_{2k-4}\mbf f_{2k-4}^x-\omega_{2k-3}\mbf f_{2k-3}^x-\nu_{2k-2}\mbf f_{2k-2}^x)/\rho_{2k-1}$
\\ \noalign{\smallskip} \quad 19:\ \  \quad  $\mbf f_{2k-1}^y= (-\xi_{2k-5}\mbf f_{2k-5}^y-\zeta_{2k-4}\mbf f_{2k-4}^y-\omega_{2k-3}\mbf f_{2k-3}^y-\nu_{2k-2}\mbf f_{2k-2}^y)/\rho_{2k-1}$ 
\\ \noalign{\smallskip} \quad 20:\ \  \quad  $\mbf f_{2k}^x= (-\xi_{2k-4}\mbf f_{2k-4}^x-\zeta_{2k-3}\mbf f_{2k-3}^x-\omega_{2k-2}\mbf f_{2k-2}^x-\nu_{2k-1}\mbf f_{2k-1}^x)/\rho_{2k}$ 
\\ \noalign{\smallskip} \quad 21:\ \  \quad $\mbf f_{2k}^y= (\mbf u_k -\xi_{2k-4}\mbf f_{2k-4}^y-\zeta_{2k-3}\mbf f_{2k-3}^y-\omega_{2k-2}\mbf f_{2k-2}^y-\nu_{2k-1}\mbf f_{2k-1}^y)/\rho_{2k}$

\\ \noalign{\smallskip} \quad 22:\ \  \quad $\mbf x_k^\rmq=\mbf x_{k-1}^\rmq+\varpi_{2k-1}\mbf f_{2k-1}^x+ \varpi_{2k}\mbf f_{2k}^x$ 
\\ \noalign{\smallskip} \quad 23:\ \  \quad $\mbf y_k^\rmq=\mbf y_{k-1}^\rmq+ \varpi_{2k-1} \mbf f_{2k-1}^y +\varpi_{2k}\mbf f_{2k}^y$
\\ \noalign{\smallskip} \quad 24:\ \  \quad $k=k+1$
\\ \noalign{\smallskip} \quad 25:\ \  {\bf end while}\\ 
\bottomrule
\end{tabular*}
\end{table}

\section{Numerical experiments}
We compare the performance of GPBiLQ, GPBiCG, GPQMR, TriCG, TriMR, and GPMR on systems generated from matrices in the SuiteSparse Matrix Collection \cite{davis2011unive}. The right-hand side of \eqref{lin} is generated so that the exact solution is the vector of ones. All algorithms stop as soon as the residual norm $\|\mbf r_k\|\leq {\tt tol}$ or $k>{\tt maxit}$, where {\tt tol} is a given tolerance and {\tt maxit} is the maximum number of iterations. To get a fair comparison, residuals are calculated explicitly at each iteration. For GPBiLQ, GPBiCG, and GPQMR, we set $\mbf f=\mbf b$ and $\mbf g=\mbf c$.

\begin{figure}[!htb]
  \centerline{\includegraphics[height=10cm]{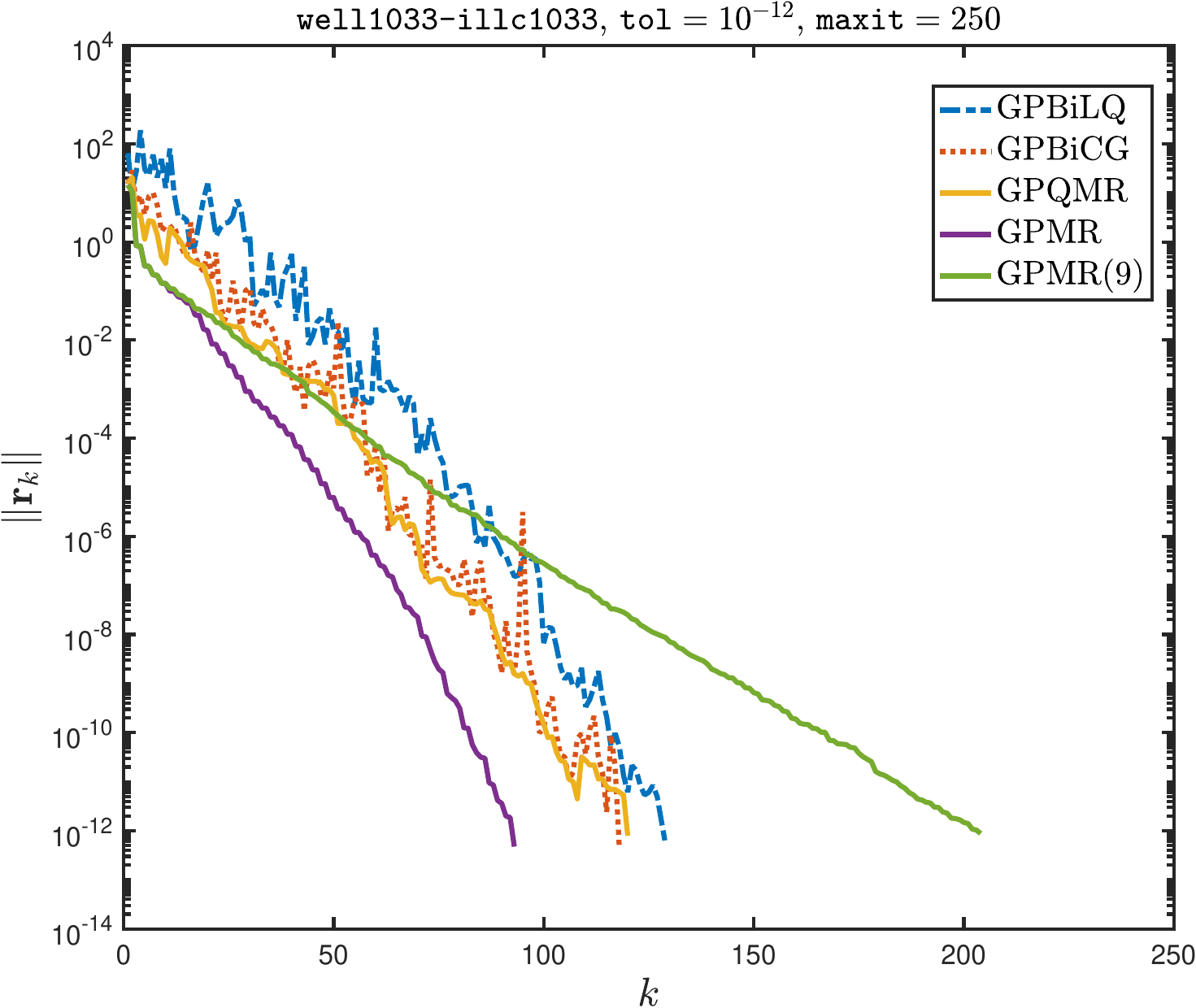}}
  \caption{Convergence curves on the unsymmetric system \eqref{lin} with matrices {\tt well1033} as $\mbf A^\top$ and {\tt illc1033} as $\mbf B$, $\lambda=1$, and $\mu=-0.1$.}\label{fig1033}
\end{figure}

\begin{figure}[!htb]
  \centerline{\includegraphics[height=10cm]{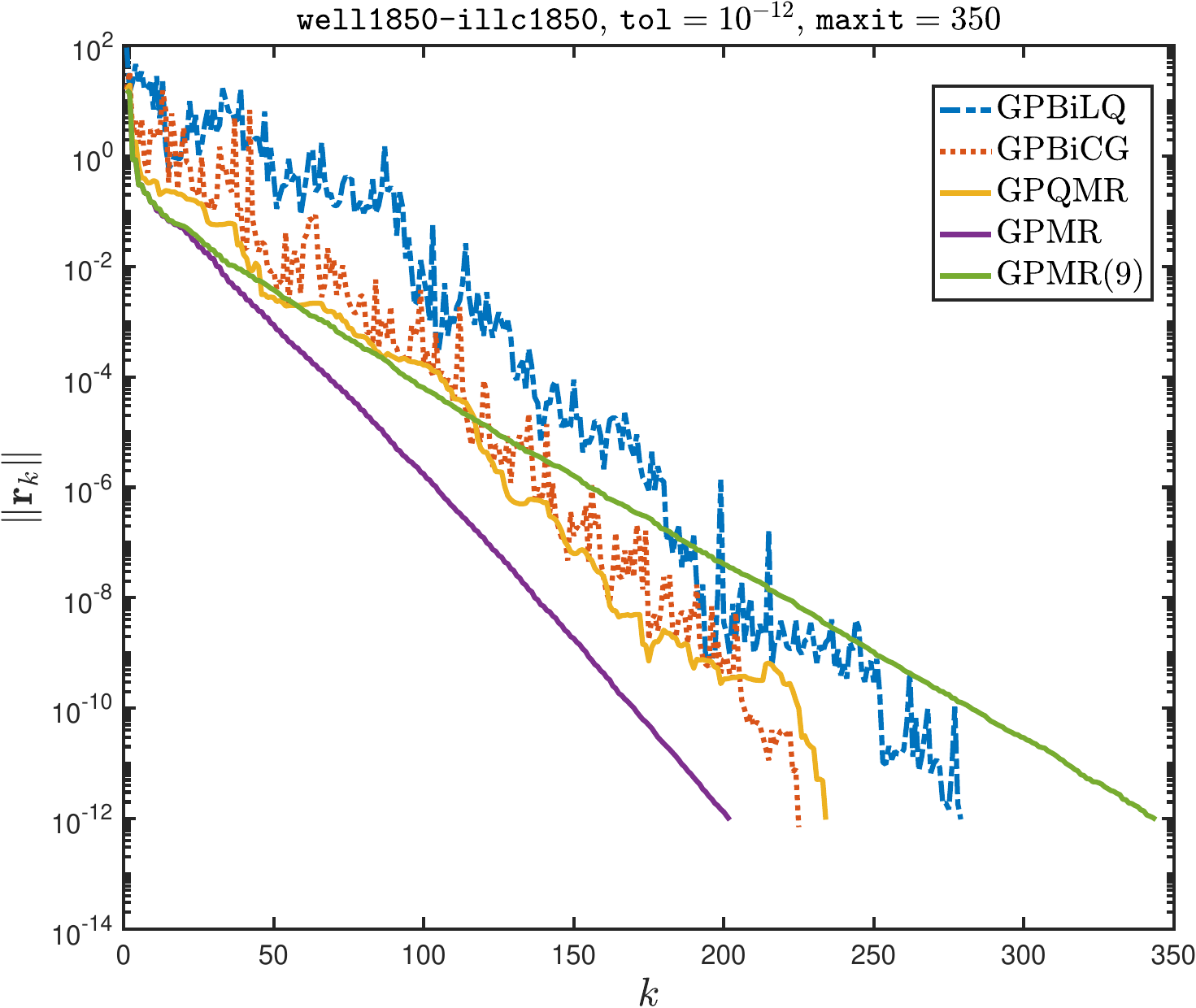}}
  \caption{Convergence curves on the unsymmetric system \eqref{lin} with matrices {\tt well1850} as $\mbf A^\top$ and {\tt illc1850} as $\mbf B$, $\lambda=1$, and $\mu=-0.05$.}\label{fig1850}
\end{figure}

We first construct two unsymmetric partitioned linear systems \eqref{lin} with different shapes of $\mbf A$ and $\mbf B$. The first system is constructed with matrices {\tt well1033} as $\mbf A^\top$ and {\tt illc1033} as $\mbf B$, $\lambda=1$, and $\mu=-0.1$. And the second one is constructed with matrices {\tt well1850} as $\mbf A^\top$ and {\tt illc1850} as $\mbf B$, $\lambda=1$, and $\mu=-0.05$. We plot the convergence curves of GPBiLQ, GPBiCG, GPQMR, GPMR, and GPMR(9) on the two unsymmetric systems in Figures \ref{fig1033} and \ref{fig1850}. We have the following observations:  
\bit
\item[(i)] GPBiLQ, GPBiCG, and GPQMR are slower than GPMR, but faster than GPMR(9);
\item[(ii)] The convergence curves of GPBiLQ, GPBiCG, and GPQMR are erratic (none of them is monotonic), although GPQMR is slightly smoother.
\eit  

\begin{figure}[!htb]
  \centerline{\includegraphics[height=10cm]{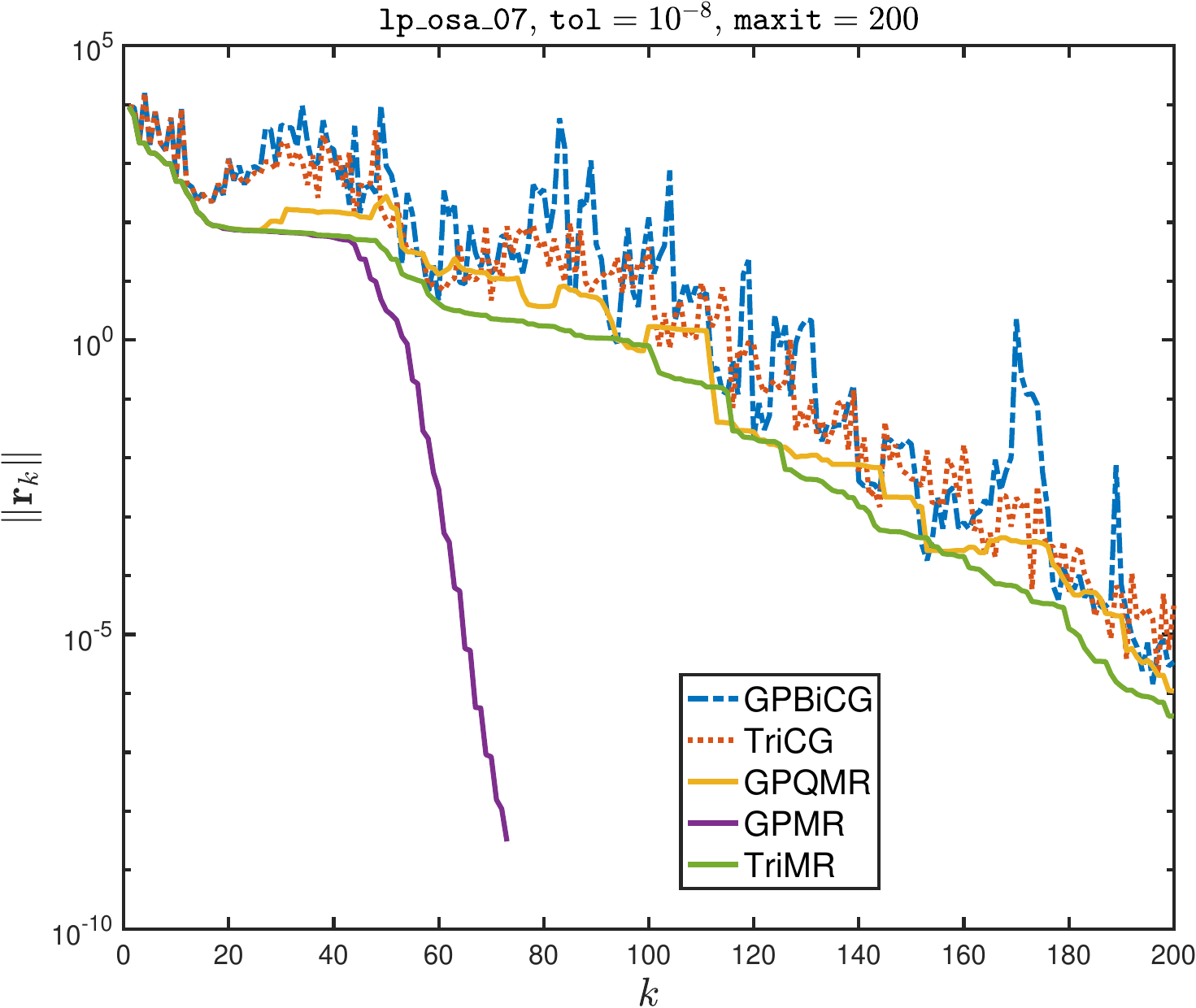}}
  \caption{Convergence curves on the symmetric system \eqref{lin} with matrix {\tt lp\_osa\_07} as $\mbf A=\mbf B^\top$, $\lambda=1$, and $\mu=-1$.}\label{figosa07}
\end{figure}

\begin{figure}[!htb]
  \centerline{\includegraphics[height=10cm]{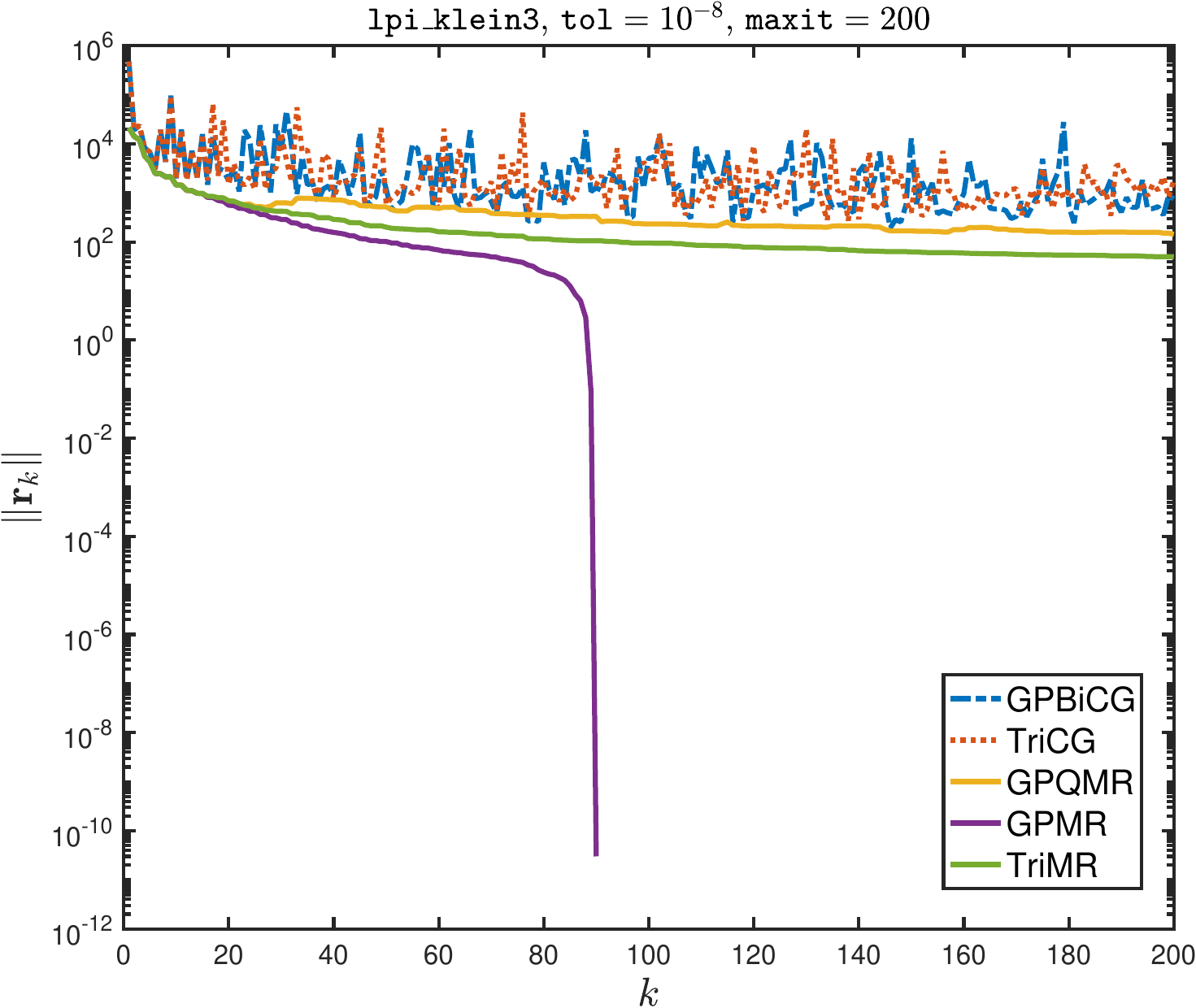}}
  \caption{Convergence curves on the symmetric system \eqref{lin} with matrix {\tt lpi\_klein3} as $\mbf A=\mbf B^\top$, $\lambda=1$, and $\mu=-1$.}\label{figklein3}
\end{figure}

By Remark \ref{relationtossy}, when applied to the system \eqref{lin} with $\mbf B=\mbf A^\top$, GPBiCG is theoretically equivalent to TriCG, and GPQMR and GPMR are theoretically equivalent to TriMR. We next evaluate the performance of these methods on two symmetric quasi-definite systems with matrices {\tt lp\_osa\_07} or {\tt lpi\_klein3} as $\mbf A$, $\lambda=1$, and $\mu=-1$. We plot the convergence curves of GPBiCG, TriCG, GPQMR, GPMR, and TriMR on the two symmetric systems in Figures \ref{figosa07} and \ref{figklein3}. We observe that GPMR is the best in terms of the number of iterations. The reason is that GPMR can be viewed as TriMR with full reorthogonalization, and the loss of orthogonality greatly affects the convergence of GPBiCG, TriCG, GPQMR, and TriMR. See \cite[section 4]{montoison2023gpmr} for more numerical comparisons between GPMR and TriMR. 

\section{Extensions and concluding remarks} Although we develop GPBiLQ and GPQMR for system \eqref{lin}, both can be generalized to solve linear systems of the form $$\bem\lambda\mbf M &\mbf A\\ \mbf B & \mu\mbf N\eem\bem\mbf x\\ \mbf y\eem=\bem\mbf b\\ \mbf c\eem,$$ where $\mbf M$ and $\mbf N$ are symmetric positive definite. We only need to introduce $\mbf M$ and $\mbf N$ in Algorithm 2 to obtain the following relations:
\begin{align*}
\mbf A \mbf U_k&=\mbf M\mbf Q_{k+1}\mbf S_{k+1,k}=\mbf M\mbf Q_k\mbf S_k+\beta_{k+1}\mbf M\mbf q_{k+1}\mbf e_k^\top,
\\ \mbf A^\top\mbf P_k &=\mbf N\mbf V_{k+1}\mbf S_{k,k+1}^\top=\mbf N\mbf V_k\mbf S_k^\top+\gamma_{k+1}\mbf N\mbf v_{k+1}\mbf e_k^\top,
\\ \mbf B \mbf Q_k &=\mbf N\mbf U_{k+1}\mbf T_{k+1,k}=\mbf N\mbf U_k\mbf T_k+\delta_{k+1}\mbf N\mbf u_{k+1}\mbf e_k^\top,
\\ \mbf B^\top\mbf V_k &=\mbf M\mbf P_{k+1}\mbf T_{k,k+1}^\top=\mbf M\mbf P_k\mbf T_k^\top+\eta_{k+1}\mbf M\mbf p_{k+1}\mbf e_k^\top.
\end{align*} Accordingly, equation \eqref{proj} becomes $$\bem\lambda\mbf M & \mbf A\\ \mbf B & \mu\mbf N\eem\bem \mbf Q_k & \mbf 0\\ \mbf 0& \mbf U_k\eem=\bem \mbf M & \mbf 0 \\ \mbf 0 & \mbf N\eem\l( \bem \mbf Q_k & \mbf 0\\ \mbf 0& \mbf U_k \eem \bem \lambda\mbf I & \mbf S_k\\ \mbf T_k & \mu\mbf I \eem + \bem \mbf q_{k+1} & \mbf 0\\ \mbf 0 & \mbf u_{k+1} \eem \bem \mbf 0 & \beta_{k+1}\mbf e_k^\top\\ \delta_{k+1}\mbf e_k^\top &  \mbf 0\eem\r),$$ 
and equation \eqref{proj2} becomes $$\bem\lambda\mbf M & \mbf A\\ \mbf B & \mu\mbf N\eem\mbf W_k=\bem \mbf M & \mbf 0 \\ \mbf 0 & \mbf N\eem\mbf W_{k+1}\mbf H_{k+1,k}.$$ Note that the biorthogonality $\mbf P_k^\top\mbf Q_k =\mbf U_k^\top\mbf V_k=\mbf I_k$ becomes the generalized biorthogonality $$\mbf P_k^\top\mbf M\mbf Q_k =\mbf U_k^\top\mbf N\mbf V_k=\mbf I_k.$$  The remaining derivation of the two methods is almost the same.

Just like GPMR can be derived by applying Block-GMRES to a special multiple right-hand-side systems (see \cite[section 5]{montoison2023gpmr}), we can derive GPBiLQ and GPQMR by Block-BiLQ and Block-QMR \cite{freund1997block}. In this paper our derivation of GPBiLQ and GPQMR follows the approach used in \cite{montoison2021tricg} and \cite{montoison2023gpmr}. 

GPBiLQ and GPQMR have one advantage over GPMR: they involve short recurrences, enabling the work per iteration and the storage requirements to remain under control even when many iterations are needed. The storage requirements of GPBiLQ and GPQMR are comparable to that of GPMR(9). For the testing examples, GPBiLQ and GPQMR are faster than GPMR(9). On the other hand, GPBiLQ and GPQMR have two disadvantages. One is that in comparison to the monotonic and often rapid convergence of GPMR as a function of iteration number, the convergence of GPBiLQ and GPQMR is slower and often erratic. The other is that except for multiplications by $\mbf A$ and $\mbf B$, GPBiLQ and GPQMR also require multiplications by their transpose. All the pros and cons that exist among GPBiLQ, GPQMR, and GPMR are similar to that among BiLQ, QMR, and GMRES. The techniques that used in CGS \cite{sonneveld1989cgs}, BiCGSTAB \cite{vandervorst1992bi}, TFQMR \cite{freund1993trans}, and BiLQR \cite{montoison2020bilq} can be used to develop variants of GPBiLQ and GPQMR. We leave them to future work.

Our MATLAB implementations of GPBiLQ, GPQMR, GPMR, TriCG, and TriMR are available at \href{https://kuidu.github.io/code.html}{https://kuidu.github.io/code.html}. 

\subsection*{Data availability} The datasets used in the numerical experiments are available from the SuiteSparse Matrix Collection (formerly known as the University of Florida Sparse Matrix Collection), \href{http://sparse.tamu.edu}{http://sparse.tamu.edu}.

\section*{Declarations}
\subsection*{Competing interests} The authors have no competing interests to declare that are relevant to the content of this article.

\appendix
\section{GPBiLQ and GPBiCG details}
\subsection{Computations for the LQ factorization of $\mbf H_k$} \label{lq}

Initialize: $$\bar\rho_1=\lambda,\quad \bar\alpha_1=\alpha_1,\quad \bar\nu_2=\theta_1,\quad \bar\rho_2=\mu,\quad \bar\omega_3=0,\quad \bar\nu_3=\beta_2,\quad \bar\zeta_4=\delta_2.$$
For $i = 1,2,\ldots,k-1$,
$$\bem\bar\rho_{2i-1} & \bar\alpha_i & 0 & \gamma_{i+1} \\ \bar\nu_{2i} & \bar\rho_{2i} & \eta_{i+1} & 0\\ \bar\omega_{2i+1} & \bar\nu_{2i+1} & \lambda & \alpha_{i+1}\\ \bar\zeta_{2i+2} & 0 & \theta_{i+1} & \mu \\ 0 & 0 & 0 & \beta_{i+2} \\ 0 & 0 & \delta_{i+2} & 0 \eem\bem c_{i,1} &&& -s_{i,1}\\ & 1&&\\ &&1&\\ s_{i,1} & & & c_{i,1}\eem=\bem\wt\rho_{2i-1} & \bar\alpha_i & 0 & 0 \\ \wt\nu_{2i} & \bar\rho_{2i} & \eta_{i+1} & t_i\\ \wt\omega_{2i+1} & \bar\nu_{2i+1} & \lambda & \wt\alpha_{i+1}\\ \wt\zeta_{2i+2} & 0 & \theta_{i+1} & \wt\rho_{2i+2} \\ \wt\xi_{2i+3} & 0 & 0 & \wt\nu_{2i+3} \\ 0 & 0 & \delta_{i+2} & 0 \eem,$$
$$\wt\rho_{2i-1}=\sqrt{\bar\rho_{2i-1}^2+\gamma_{i+1}^2},\quad c_{i,1}=\bar\rho_{2i-1}/\wt\rho_{2i-1},\quad s_{i,1}=\gamma_{i+1}/\wt\rho_{2i-1},$$
$$\wt\nu_{2i}= c_{i,1}\bar\nu_{2i},\quad t_i=-s_{i,1}\bar\nu_{2i},$$ $$ \wt\omega_{2i+1}=c_{i,1}\bar\omega_{2i+1}+s_{i,1}\alpha_{i+1},\quad \wt\alpha_{i+1}=-s_{i,1}\bar\omega_{2i+1}+c_{i,1}\alpha_{i+1},$$ 
$$ \wt\zeta_{2i+2}=c_{i,1}\bar\zeta_{2i+2}+s_{i,1}\mu,\quad \wt\rho_{2i+2}=-s_{i,1}\bar\zeta_{2i+2}+c_{i,1}\mu,$$ $$ \wt\xi_{2i+3}=s_{i,1}\beta_{i+2},\quad  \wt\nu_{2i+3}=c_{i,1}\beta_{i+2},$$
$$\bem\wt\rho_{2i-1} & \bar\alpha_i & 0 & 0 \\ \wt\nu_{2i} & \bar\rho_{2i} & \eta_{i+1} & t_i\\ \wt\omega_{2i+1} & \bar\nu_{2i+1} & \lambda & \wt\alpha_{i+1}\\ \wt\zeta_{2i+2} & 0 & \theta_{i+1} & \wt\rho_{2i+2} \\ \wt\xi_{2i+3} & 0 & 0 & \wt\nu_{2i+3} \\ 0 & 0 & \delta_{i+2} & 0 \eem\bem c_{i,2}&-s_{i,2}&& \\ s_{i,2}&c_{i,2}&&\\ &&1&\\ &&&1 \eem=\bem\rho_{2i-1} & 0 & 0 & 0 \\ \nu_{2i} & \wh\rho_{2i} & \eta_{i+1} & t_i\\ \omega_{2i+1} & \wh\nu_{2i+1} & \lambda & \wt\alpha_{i+1}\\ \zeta_{2i+2} & \wh\omega_{2i+2} & \theta_{i+1} & \wt\rho_{2i+2} \\ \xi_{2i+3} & \wh\zeta_{2i+3} & 0 & \wt\nu_{2i+3} \\ 0 & 0 & \delta_{i+2} & 0 \eem,$$
$$\rho_{2i-1}=\sqrt{\wt\rho_{2i-1}^2+\bar\alpha_i^2}, \quad c_{i,2}=\wt\rho_{2i-1}/\rho_{2i-1},\quad s_{i,2}=\bar\alpha_i/\rho_{2i-1},$$
$$\nu_{2i}=c_{i,2}\wt\nu_{2i}+s_{i,2}\bar\rho_{2i},\quad  \wh\rho_{2i}=-s_{i,2}\wt\nu_{2i}+c_{i,2}\bar\rho_{2i}, $$ 
$$\omega_{2i+1}=c_{i,2}\wt\omega_{2i+1}+s_{i,2}\bar\nu_{2i+1}, \quad \wh\nu_{2i+1}=-s_{i,2}\wt\omega_{2i+1}+c_{i,2}\bar\nu_{2i+1},$$
$$\zeta_{2i+2} =c_{i,2}\wt\zeta_{2i+2}, \quad \wh\omega_{2i+2}=-s_{i,2}\wt\zeta_{2i+2},$$ $$\xi_{2i+3}=c_{i,2}\wt\xi_{2i+3},\quad \wh\zeta_{2i+3}=-s_{i,2}\wt\xi_{2i+3},$$
$$\bem\rho_{2i-1} & 0 & 0 & 0 \\ \nu_{2i} & \wh\rho_{2i} & \eta_{i+1} & t_i\\ \omega_{2i+1} & \wh\nu_{2i+1} & \lambda & \wt\alpha_{i+1}\\ \zeta_{2i+2} & \wh\omega_{2i+2} & \theta_{i+1} & \wt\rho_{2i+2} \\ \xi_{2i+3} & \wh\zeta_{2i+3} & 0 & \wt\nu_{2i+3} \\ 0 & 0 & \delta_{i+2} & 0 \eem\bem 1&&& \\ &c_{i,3}&&-s_{i,3}\\ &&1&\\ &s_{i,3}&&c_{i,3} \eem=\bem\rho_{2i-1} & 0 & 0 & 0 \\ \nu_{2i} & \wc\rho_{2i} & \eta_{i+1} & 0\\ \omega_{2i+1} & \wc\nu_{2i+1} & \lambda & \bar\alpha_{i+1}\\ \zeta_{2i+2} & \wc\omega_{2i+2} & \theta_{i+1} & \bar\rho_{2i+2} \\ \xi_{2i+3} & \wc\zeta_{2i+3} & 0 & \bar\nu_{2i+3} \\ 0 & 0 & \delta_{i+2} & 0 \eem,$$
$$\wc\rho_{2i}=\sqrt{\wh\rho_{2i}^2+t_i^2},\quad c_{i,3}=\wh\rho_{2i}/\wc\rho_{2i},\quad s_{i,3}=t_i/\wc\rho_{2i},$$
$$\wc \nu_{2i+1}=c_{i,3}\wh\nu_{2i+1}+s_{i,3}\wt\alpha_{i+1},\quad \bar\alpha_{i+1}=-s_{i,3}\wh\nu_{2i+1}+c_{i,3}\wt\alpha_{i+1},$$ $$ \wc\omega_{2i+2}=c_{i,3}\wh\omega_{2i+2}+s_{i,3}\wt\rho_{2i+2},\quad \bar\rho_{2i+2}=-s_{i,3}\wh\omega_{2i+2}+c_{i,3}\wt\rho_{2i+2},$$
$$\wc\zeta_{2i+3}=c_{i,3}\wh\zeta_{2i+3}+s_{i,3}\wt\nu_{2i+3},\quad \bar\nu_{2i+3}=-s_{i,3}\wh\zeta_{2i+3}+c_{i,3}\wt\nu_{2i+3}, $$
$$\bem\rho_{2i-1} & 0 & 0 & 0 \\ \nu_{2i} & \wc\rho_{2i} & \eta_{i+1} & 0\\ \omega_{2i+1} & \wc\nu_{2i+1} & \lambda & \bar\alpha_{i+1}\\ \zeta_{2i+2} & \wc\omega_{2i+2} & \theta_{i+1} & \bar\rho_{2i+2} \\ \xi_{2i+3} & \wc\zeta_{2i+3} & 0 & \bar\nu_{2i+3} \\ 0 & 0 & \delta_{i+2} & 0 \eem\bem 1&&&\\ &c_{i,4}& -s_{i,4}& \\ & s_{i,4}& c_{i,4}&\\ &&&1\eem =\bem\rho_{2i-1} & 0 & 0 & 0 \\ \nu_{2i} & \rho_{2i} & 0 & 0\\ \omega_{2i+1} &\nu_{2i+1} & \bar\rho_{2i+1} & \bar\alpha_{i+1}\\ \zeta_{2i+2} & \omega_{2i+2} & \bar\nu_{2i+2} & \bar\rho_{2i+2} \\ \xi_{2i+3} & \zeta_{2i+3} & \bar\omega_{2i+3} & \bar\nu_{2i+3} \\ 0 & \xi_{2i+4} & \bar\zeta_{2i+4} & 0 \eem,$$
$$\rho_{2i}=\sqrt{\wc\rho_{2i}^2+\eta_{i+1}^2},\quad c_{i,4}=\wc\rho_{2i}/\rho_{2i},\quad s_{i,4}=\eta_{i+1}/\rho_{2i},$$
$$\nu_{2i+1}=c_{i,4}\wc\nu_{2i+1}+s_{i,4}\lambda,\quad \bar\rho_{2i+1}=-s_{i,4}\wc\nu_{2i+1}+c_{i,4}\lambda, $$
$$\omega_{2i+2}=c_{i,4}\wc\omega_{2i+2}+s_{i,4}\theta_{i+1},\quad \bar\nu_{2i+2}=-s_{i,4}\wc\omega_{2i+2}+c_{i,4}\theta_{i+1},$$
$$\zeta_{2i+3}=c_{i,4}\wc\zeta_{2i+3},\quad \bar\omega_{2i+3}=-s_{i,4}\wc\zeta_{2i+3},$$
$$\xi_{2i+4}=s_{i,4}\delta_{i+2},\quad \bar\zeta_{2i+4}=c_{i,4}\delta_{i+2}.$$
The structure of $\wt{\mbf G}_{2k-1,2k}$ is 
$$\wt{\mbf G}_{2k-1,2k}=\bem\mbf I_{2k-2} &&\\ &c_k&-s_k\\ &s_k& c_k\eem.$$ We have
$$\bem \bar\rho_{2k-1} & \bar\alpha_k \\  \bar\nu_{2k} & \bar\rho_{2k} \eem \bem c_k&-s_k\\ s_k& c_k\eem=\bem  \ddot\rho_{2k-1} & 0 \\ \ddot\nu_{2k} & \ddot\rho_{2k} \eem,$$
$$\ddot\rho_{2k-1}=\sqrt{\bar\rho_{2k-1}^2+\bar\alpha_k^2},\quad c_k=\bar\rho_{2k-1}/\ddot\rho_{2k-1},\quad s_k=\bar\alpha_k/\ddot\rho_{2k-1},$$ 
$$\ddot\nu_{2k}=c_k\bar\nu_{2k}+s_k\bar\rho_{2k},\quad \ddot\rho_{2k}=-s_k\bar\nu_{2k}+c_k\bar\rho_{2k}.$$

\subsection{Computations for $\varpi_1$, $\ldots$, $\varpi_{2k-2}$, $\wt\varpi_{2k-1}$, and $\wt\varpi_{2k}$}\label{wtt}
We have$$\bem \varpi_1 & \ldots &\varpi_{2k-2} & \wt\varpi_{2k-1} & \wt\varpi_{2k}\eem^\top =\wt{\mbf L}_k^{-1}(\beta_1\mbf e_1+\delta_1\mbf e_2),$$
$$\varpi_1=\beta_1/\rho_1,\quad \varpi_2=(\delta_1-\nu_2\varpi_1)/\rho_2,\quad \varpi_3=-(\omega_3\varpi_1+\nu_3\varpi_2)/\rho_3,\quad\varpi_4=-(\zeta_4\varpi_1+\omega_4\varpi_2+\nu_4\varpi_3)/\rho_4,$$
$$\varpi_5=-(\xi_5\varpi_1+\zeta_5\varpi_2+\omega_5\varpi_3+\nu_5\varpi_4)/\rho_5,\qquad \cdots,$$
$$\varpi_{2k-3}=-(\xi_{2k-3}\varpi_{2k-7}+\zeta_{2k-3}\varpi_{2k-6}+\omega_{2k-3}\varpi_{2k-5}+\nu_{2k-3}\varpi_{2k-4})/\rho_{2k-3}\quad (k\geq 4),$$
$$\varpi_{2k-2}=-(\xi_{2k-2}\varpi_{2k-6}+\zeta_{2k-2}\varpi_{2k-5}+\omega_{2k-2}\varpi_{2k-4}+\nu_{2k-2}\varpi_{2k-3})/\rho_{2k-2}\quad (k\geq 4),$$
$$\wt\varpi_{2k-1}=-(\xi_{2k-1}\varpi_{2k-5}+\zeta_{2k-1}\varpi_{2k-4}+\omega_{2k-1}\varpi_{2k-3}+\nu_{2k-1}\varpi_{2k-2})/\ddot\rho_{2k-1},$$ 
$$\wt\varpi_{2k}=-(\xi_{2k}\varpi_{2k-4}+\zeta_{2k}\varpi_{2k-3}+\omega_{2k}\varpi_{2k-2}+\ddot\nu_{2k}\wt\varpi_{2k-1})/\ddot\rho_{2k}.$$

\subsection{Computations for $\vartheta_k$, $\varrho_k$, $\chi_k$, and $\varsigma_k$}\label{trcs}
Because $\mbf G_k=\mbf G_{1,4}\mbf G_{3,6}\cdots\mbf G_{2k-5,2k-2}\mbf G_{2k-3,2k}$ and $\mbf z_k^\rml=\mbf G_k\bem\mbf t_{k-1}\\ \mbf 0\eem$, we have \begin{align*}\bem\mbf e_{2k-3}^\top\\ \mbf e_{2k-2}^\top \eem\mbf z_k^\rml &=\bem\mbf e_{2k-3}^\top\\ \mbf e_{2k-2}^\top \eem\mbf G_{1,4}\mbf G_{3,6}\cdots\mbf G_{2k-5,2k-2}\mbf G_{2k-3,2k}\bem\mbf t_{k-1}\\ \mbf 0\eem \\ &= \bem\mbf e_{2k-3}^\top\\ \mbf e_{2k-2}^\top \eem \mbf G_{2k-5,2k-2}\mbf G_{2k-3,2k}\bem\mbf t_{k-1}\\ \mbf 0\eem \\ & = \bem\mbf e_{2k-3}^\top\\ \mbf e_{2k-2}^\top \eem \mbf G_{2k-5,2k-2}\mbf G_{2k-3,2k}\bem \mbf e_{2k-5} & \mbf e_{2k-4} & \cdots & \mbf e_{2k} \eem \bem \varpi_{2k-5}\\ \varpi_{2k-4} \\ \varpi_{2k-3} \\ \varpi_{2k-2} \\   0 \\ 0 \eem \end{align*}
and 
\begin{align*}
	\bem\mbf e_{2k-1}^\top\\ \mbf e_{2k}^\top \eem\mbf z_k^\rml &=\bem\mbf e_{2k-1}^\top\\ \mbf e_{2k}^\top \eem\mbf G_{1,4}\mbf G_{3,6}\cdots\mbf G_{2k-5,2k-2}\mbf G_{2k-3,2k}\bem\mbf t_{k-1}\\ \mbf 0\eem \\ &= \bem\mbf e_{2k-1}^\top\\ \mbf e_{2k}^\top \eem \mbf G_{2k-3,2k}\bem\mbf t_{k-1}\\ \mbf 0\eem\\ & = \bem\mbf e_{2k-1}^\top\\ \mbf e_{2k}^\top \eem \mbf G_{2k-3,2k}\bem \mbf e_{2k-3} & \mbf e_{2k-2} & \mbf e_{2k-1} & \mbf e_{2k} \eem\bem  \varpi_{2k-3} \\ \varpi_{2k-2} \\   0 \\ 0 \eem.	
\end{align*} 
The scalars $\vartheta_k$, $\varrho_k$, $\chi_k$, and $\varsigma_k$ can be computed via $$\bem \vartheta_k \\ \varrho_k \eem= \l(\mbf E_{k,k-1}\bem\mbf e_{2k-3}^\top\\ \mbf e_{2k-2}^\top \eem + \mbf E_{k,k}\bem\mbf e_{2k-1}^\top\\ \mbf e_{2k}^\top \eem\r) \mbf z_k^\rml,\qquad\bem\chi_k\\ \varsigma_k\eem= \mbf E_{k+1,k}\bem\mbf e_{2k-1}^\top\\ \mbf e_{2k}^\top \eem\mbf z_k^\rml.$$

\subsection{Computations for $\wt\chi_k$ and $\wt\varsigma_k$}\label{tcs}
Because $\wt{\mbf Q}_k^\top=\mbf G_{1,4}\mbf G_{3,6}\cdots\mbf G_{2k-5,2k-2}\mbf G_{2k-3,2k}\wt{\mbf G}_{2k-1,2k}$ and $\mbf z_k^\rmc=\wt{\mbf Q}_k^\top\wt{\mbf t}_k$, we have
\begin{align*}
\bem\mbf e_{2k-1}^\top\\ \mbf e_{2k}^\top \eem\mbf z_k^\rmc & =\bem\mbf e_{2k-1}^\top\\ \mbf e_{2k}^\top \eem\mbf G_{1,4}\mbf G_{3,6}\cdots\mbf G_{2k-5,2k-2}\mbf G_{2k-3,2k}\wt{\mbf G}_{2k-1,2k}\wt{\mbf t}_k\\ &=\bem\mbf e_{2k-1}^\top\\ \mbf e_{2k}^\top \eem \mbf G_{2k-3,2k}\wt{\mbf G}_{2k-1,2k}\wt{\mbf t}_k\\ & = \bem\mbf e_{2k-1}^\top\\ \mbf e_{2k}^\top \eem \mbf G_{2k-3,2k}\wt{\mbf G}_{2k-1,2k}\bem \mbf e_{2k-3} & \mbf e_{2k-2} & \mbf e_{2k-1} & \mbf e_{2k} \eem\bem  \varpi_{2k-3} \\ \varpi_{2k-2} \\  \wt\varpi_{2k-1} \\ \wt\varpi_{2k}\eem.
\end{align*}
The scalars $\wt\chi_k$ and $\wt\varsigma_k$ can be computed via $$\bem\wt\chi_k\\ \wt\varsigma_k\eem= \mbf E_{k+1,k}\bem\mbf e_{2k-1}^\top\\ \mbf e_{2k}^\top \eem\mbf z_k^\rmc.$$
 
\section{GPQMR details}
\subsection{Computations for the QR factorization of $\mbf H_{k+1,k}$}\label{qr}
Initialize: $$\bar\rho_1=\lambda,\quad \bar\theta_1=\theta_1,\quad \bar\nu_1=\alpha_1,\quad \bar\rho_2=\mu,\quad \bar\omega_1=0,\quad \bar\nu_2=\eta_2,\quad \bar\zeta_1=\gamma_2.$$
For $i = 1,2,\ldots,k$,
{\small
$$\bem c_{i,1} &&& s_{i,1}\\ & 1&&\\ &&1&\\ -s_{i,1} & & & c_{i,1}\eem\bem \bar\rho_{2i-1} & \bar\nu_{2i-1} & \bar\omega_{2i-1} & \bar\zeta_{2i-1} & 0& 0\\ \bar\theta_i & \bar\rho_{2i} & \bar\nu_{2i} & 0 & 0 & 0 \\ 0 & \beta_{i+1} & \lambda & \alpha_{i+1} & 0 & \gamma_{i+2}\\ \delta_{i+1} & 0 & \theta_{i+1} & \mu & \eta_{i+2} & 0\eem=\bem \wt\rho_{2i-1} & \wt\nu_{2i-1} & \wt\omega_{2i-1} & \wt\zeta_{2i-1} & \wt\xi_{2i-1} & 0\\ \bar\theta_i & \bar\rho_{2i} & \bar\nu_{2i} & 0 & 0 & 0 \\ 0 & \beta_{i+1} & \lambda & \alpha_{i+1} & 0 & \gamma_{i+2}\\ 0 & t_i & \wt\theta_{i+1} & \wt\rho_{2i+2} & \wt\nu_{2i+2} & 0\eem,$$}
$$\wt\rho_{2i-1}=\sqrt{\bar\rho_{2i-1}^2+\delta_{i+1}^2},\quad c_{i,1}=\bar\rho_{2i-1}/\wt\rho_{2i-1},\quad s_{i,1}=\delta_{i+1}/\wt\rho_{2i-1},$$
$$\wt\nu_{2i-1}= c_{i,1}\bar\nu_{2i-1},\quad t_i=-s_{i,1}\bar\nu_{2i-1},$$ $$ \wt\omega_{2i-1}=c_{i,1}\bar\omega_{2i-1}+s_{i,1}\theta_{i+1},\quad \wt\theta_{i+1}=-s_{i,1}\bar\omega_{2i-1}+c_{i,1}\theta_{i+1},$$ 
$$ \wt\zeta_{2i-1}=c_{i,1}\bar\zeta_{2i-1}+s_{i,1}\mu,\quad \wt\rho_{2i+2}=-s_{i,1}\bar\zeta_{2i-1}+c_{i,1}\mu,$$ $$ \wt\xi_{2i-1}=s_{i,1}\eta_{i+2},\quad  \wt\nu_{2i+2}=c_{i,1}\eta_{i+2},$$
{\small
$$\bem c_{i,2}&s_{i,2}&& \\ -s_{i,2}&c_{i,2}&&\\ &&1&\\ &&&1 \eem\bem \wt\rho_{2i-1} & \wt\nu_{2i-1} & \wt\omega_{2i-1} & \wt\zeta_{2i-1} & \wt\xi_{2i-1} & 0\\ \bar\theta_i & \bar\rho_{2i} & \bar\nu_{2i} & 0 & 0 & 0 \\ 0 & \beta_{i+1} & \lambda & \alpha_{i+1} & 0 & \gamma_{i+2}\\ 0 & t_i & \wt\theta_{i+1} & \wt\rho_{2i+2} & \wt\nu_{2i+2} & 0\eem=\bem \rho_{2i-1} & \nu_{2i-1} & \omega_{2i-1} & \zeta_{2i-1} & \xi_{2i-1} & 0\\ 0 & \wh\rho_{2i} & \wh\nu_{2i} & \wh\omega_{2i} & \wh\zeta_{2i} & 0 \\ 0 & \beta_{i+1} & \lambda & \alpha_{i+1} & 0 & \gamma_{i+2}\\ 0 & t_i & \wt\theta_{i+1} & \wt\rho_{2i+2} & \wt\nu_{2i+2} & 0\eem,$$}
$$\rho_{2i-1}=\sqrt{\wt\rho_{2i-1}^2+\bar\theta_i^2}, \quad c_{i,2}=\wt\rho_{2i-1}/\rho_{2i-1},\quad s_{i,2}=\bar\theta_i/\rho_{2i-1},$$
$$\nu_{2i-1}=c_{i,2}\wt\nu_{2i-1}+s_{i,2}\bar\rho_{2i},\quad  \wh\rho_{2i}=-s_{i,2}\wt\nu_{2i-1}+c_{i,2}\bar\rho_{2i}, $$ 
$$\omega_{2i-1}=c_{i,2}\wt\omega_{2i-1}+s_{i,2}\bar\nu_{2i}, \quad \wh\nu_{2i}=-s_{i,2}\wt\omega_{2i-1}+c_{i,2}\bar\nu_{2i},$$
$$\zeta_{2i-1} =c_{i,2}\wt\zeta_{2i-1}, \quad \wh\omega_{2i}=-s_{i,2}\wt\zeta_{2i-1},$$ $$\xi_{2i-1}=c_{i,2}\wt\xi_{2i-1},\quad \wh\zeta_{2i}=-s_{i,2}\wt\xi_{2i-1},$$
{\small
$$\bem 1&&& \\ &c_{i,3}&&s_{i,3}\\ &&1&\\ &-s_{i,3}&&c_{i,3} \eem\bem \rho_{2i-1} & \nu_{2i-1} & \omega_{2i-1} & \zeta_{2i-1} & \xi_{2i-1} & 0\\ 0 & \wh\rho_{2i} & \wh\nu_{2i} & \wh\omega_{2i} & \wh\zeta_{2i} & 0 \\ 0 & \beta_{i+1} & \lambda & \alpha_{i+1} & 0 & \gamma_{i+2}\\ 0 & t_i & \wt\theta_{i+1} & \wt\rho_{2i+2} & \wt\nu_{2i+2} & 0\eem=\bem \rho_{2i-1} & \nu_{2i-1} & \omega_{2i-1} & \zeta_{2i-1} & \xi_{2i-1} & 0\\ 0 & \wc\rho_{2i} & \wc\nu_{2i} & \wc\omega_{2i} & \wc\zeta_{2i} & 0 \\ 0 & \beta_{i+1} & \lambda & \alpha_{i+1} & 0 & \gamma_{i+2}\\ 0 & 0 & \bar\theta_{i+1} & \bar\rho_{2i+2} & \bar\nu_{2i+2} & 0\eem,$$}
$$\wc\rho_{2i}=\sqrt{\wh\rho_{2i}^2+t_i^2},\quad c_{i,3}=\wh\rho_{2i}/\wc\rho_{2i},\quad s_{i,3}=t_i/\wc\rho_{2i},$$
$$\wc \nu_{2i}=c_{i,3}\wh\nu_{2i}+s_{i,3}\wt\theta_{i+1},\quad \bar\theta_{i+1}=-s_{i,3}\wh\nu_{2i}+c_{i,3}\wt\theta_{i+1},$$ $$ \wc\omega_{2i}=c_{i,3}\wh\omega_{2i}+s_{i,3}\wt\rho_{2i+2},\quad \bar\rho_{2i+2}=-s_{i,3}\wh\omega_{2i}+c_{i,3}\wt\rho_{2i+2},$$
$$\wc\zeta_{2i}=c_{i,3}\wh\zeta_{2i}+s_{i,3}\wt\nu_{2i+2},\quad \bar\nu_{2i+2}=-s_{i,3}\wh\zeta_{2i}+c_{i,3}\wt\nu_{2i+2},$$
{\small
$$\bem 1&&&\\ &c_{i,4}& s_{i,4}& \\ & -s_{i,4}& c_{i,4}&\\ &&&1\eem\bem \rho_{2i-1} & \nu_{2i-1} & \omega_{2i-1} & \zeta_{2i-1} & \xi_{2i-1} & 0\\ 0 & \wc\rho_{2i} & \wc\nu_{2i} & \wc\omega_{2i} & \wc\zeta_{2i} & 0 \\ 0 & \beta_{i+1} & \lambda & \alpha_{i+1} & 0 & \gamma_{i+2}\\ 0 & 0 & \bar\theta_{i+1} & \bar\rho_{2i+2} & \bar\nu_{2i+2} & 0\eem =\bem \rho_{2i-1} & \nu_{2i-1} & \omega_{2i-1} & \zeta_{2i-1} & \xi_{2i-1} & 0\\ 0 & \rho_{2i} & \nu_{2i} & \omega_{2i} & \zeta_{2i} & \xi_{2i} \\ 0 & 0 & \bar\rho_{2i+1} & \bar\nu_{2i+1} & \bar\omega_{2i+1} & \bar\zeta_{2i+1}\\ 0 & 0 & \bar\theta_{i+1} & \bar\rho_{2i+2} & \bar\nu_{2i+2} & 0\eem,$$}
$$\rho_{2i}=\sqrt{\wc\rho_{2i}^2+\beta_{i+1}^2},\quad c_{i,4}=\wc\rho_{2i}/\rho_{2i},\quad s_{i,4}=\beta_{i+1}/\rho_{2i},$$
$$\nu_{2i}=c_{i,4}\wc\nu_{2i}+s_{i,4}\lambda,\quad \bar\rho_{2i+1}=-s_{i,4}\wc\nu_{2i}+c_{i,4}\lambda, $$
$$\omega_{2i}=c_{i,4}\wc\omega_{2i}+s_{i,4}\alpha_{i+1},\quad \bar\nu_{2i+1}=-s_{i,4}\wc\omega_{2i}+c_{i,4}\alpha_{i+1},$$
$$\zeta_{2i}=c_{i,4}\wc\zeta_{2i},\quad \bar\omega_{2i+1}=-s_{i,4}\wc\zeta_{2i},$$
$$\xi_{2i}=s_{i,4}\gamma_{i+2},\quad \bar\zeta_{2i+1}=c_{i,4}\gamma_{i+2}.$$

\subsection{Computations for $\varpi_1$, $\ldots$, $\varpi_{2k}$, $\bar\varpi_{2k+1}$, and $\bar\varpi_{2k+2}$}\label{bart} 
We have $$\bem \varpi_1 & \ldots &\varpi_{2k} & \bar\varpi_{2k+1} & \bar\varpi_{2k+2}\eem^\top =\wh{\mbf Q}_k^\top(\beta_1\mbf e_1+\delta_1\mbf e_2)=\wh{\mbf G}_{2k-1,2k+2}\cdots\wh{\mbf G}_{3,6}\wh{\mbf G}_{1,4}(\beta_1\mbf e_1+\delta_1\mbf e_2).$$
Initialize: $\bar\varpi_1=\beta_1$, $ \bar\varpi_2=\delta_1$.
For $i=1,2,\ldots,k$,
$$\bem c_{i,1} &&& s_{i,1}\\ & 1&&\\ &&1&\\ -s_{i,1} & & & c_{i,1}\eem\bem \bar\varpi_{2i-1} \\ \bar\varpi_{2i} \\ 0 \\ 0\eem= \bem \wt\varpi_{2i-1} \\ \bar\varpi_{2i} \\ 0 \\ \wt\varpi_{2i+2}\eem,$$
$$\wt\varpi_{2i-1}=c_{i,1}\bar\varpi_{2i-1},\quad  \wt\varpi_{2i+2}=-s_{i,1}\bar\varpi_{2i-1},$$
$$
\bem c_{i,2}&s_{i,2}&& \\ -s_{i,2}&c_{i,2}&&\\ &&1&\\ &&&1 \eem \bem \wt\varpi_{2i-1} \\ \bar\varpi_{2i} \\ 0 \\ \wt\varpi_{2i+2}\eem= \bem \varpi_{2i-1} \\ \wh\varpi_{2i} \\ 0 \\ \wt\varpi_{2i+2}\eem,
$$
$$\varpi_{2i-1}=c_{i,2}\wt\varpi_{2i-1}+ s_{i,2}\bar\varpi_{2i},\quad \wh\varpi_{2i}=-s_{i,2}\wt\varpi_{2i-1}+ c_{i,2}\bar\varpi_{2i},$$
$$\bem 1&&& \\ &c_{i,3}&&s_{i,3}\\ &&1&\\ &-s_{i,3}&&c_{i,3} \eem\bem \varpi_{2i-1} \\ \wh\varpi_{2i} \\ 0 \\ \wt\varpi_{2i+2}\eem=\bem \varpi_{2i-1} \\ \wc\varpi_{2i} \\ 0 \\ \bar\varpi_{2i+2}\eem,$$
$$\wc\varpi_{2i}=c_{i,3}\wh\varpi_{2i}+s_{i,3}\wt\varpi_{2i+2},\quad \bar\varpi_{2i+2}=-s_{i,3}\wh\varpi_{2i}+c_{i,3}\wt\varpi_{2i+2},$$
$$\bem 1&&&\\ &c_{i,4}& s_{i,4}& \\ & -s_{i,4}& c_{i,4}&\\ &&&1 \eem \bem \varpi_{2i-1} \\ \wc\varpi_{2i} \\ 0 \\ \bar\varpi_{2i+2}\eem=\bem \varpi_{2i-1} \\ \varpi_{2i} \\ \bar\varpi_{2i+1} \\ \bar\varpi_{2i+2}\eem,$$
$$\varpi_{2i}=c_{i,4}\wc\varpi_{2i},\quad \bar\varpi_{2i+1}=-s_{i,4}\wc\varpi_{2i}.$$

{\small 

}
\end{document}